\documentclass[12pt]{amsart}
\usepackage{amssymb}
\usepackage{epsf, epsfig}
\usepackage{multirow}
\usepackage{colordvi}

\numberwithin{equation}{section}

\newtheorem{thm}[equation]{Theorem}
\newtheorem{prop}[equation]{Proposition}
\newtheorem{lemma}[equation]{Lemma}
\newtheorem{cor}[equation]{Corollary}

\newtheorem{example}[equation]{Example}
\newtheorem{remark}[equation]{Remark}
\newtheorem{definition}[equation]{Definition}

\newenvironment{defn}{\begin{definition}\rm}{\end{definition}}
\newenvironment{ex}{\begin{example}\rm}{\end{example}}
\newenvironment{rem}{\begin{remark}\rm}{\end{remark}}
\newcounter{FNC}[page]
\def\fauxfootnote#1{{\addtocounter{FNC}{2}$^\fnsymbol{FNC}$%
     \let\thefootnote\relax\footnotetext{$^\fnsymbol{FNC}$#1}}}

\newcommand{\ch}{\mbox{\rm char}}
\newcommand{\sign}{\mbox{\,\rm sign\,}}
\newcommand{\ini}{\mbox{\rm in}}
\newcommand{\calA}{\mathcal{A}}

\newcommand{\calJ}{\mathcal{J}}
\newcommand{\calC}{\mathcal{C}}

\newcommand{\C}{\mathbb{C}}
\newcommand{\RP}{\mathbb{RP}}
\newcommand{\Q}{\mathbb{Q}}
\newcommand{\R}{\mathbb{R}}
\newcommand{\Z}{\mathbb{Z}}
\renewcommand{\P}{\mathbb{P}}

\copyrightinfo{}{}
\title[Lower Bounds for Real Polynomial Systems]{Lower Bounds for Real Solutions
  to Sparse Polynomial Systems}  

\author{Evgenia Soprunova}
\address{Department of Mathematics\\
        University of Massachusetts\\
        Amherst, MA, 01003\\
        USA}
\email{esoprun@math.umass.edu}
\urladdr{http://www.math.umass.edu/\~{}esoprun}

\author{Frank Sottile}
\address{Department of Mathematics\\
         Texas A\&M University\\
         College Station\\
         TX \ 77843\\
         USA}
\email{sottile@math.tamu.edu}
\urladdr{http://www.math.tamu.edu/\~{}sottile}

\thanks{Work of Sottile is supported by the Clay Mathematical Institute and this
  was completed while in residence at the MSRI}
\thanks{Both authors were supported in part by NSF CAREER grant DMS-0134860.}

\begin{document}

\begin{abstract}
 We show how to construct sparse polynomial systems that have 
 non-trivial lower bounds on their numbers of real solutions.
 These are unmixed systems associated to certain polytopes.
 For the order polytope of a poset $P$ this lower bound is the sign-imbalance
 of $P$ and it holds if all maximal chains of $P$ have length of the same parity.
 This theory also gives lower bounds in the real Schubert calculus
 through the sagbi
 degeneration of the Grassmannian to a toric variety, and thus recovers a result
 of Eremenko and Gabrielov.
 
\end{abstract}

\maketitle

\noindent{\it Dedicated to Richard P.~Stanley on the occasion of his
              60th Birthday.}
\section*{Introduction}

A fundamental problem in real algebraic geometry is to understand the real
solutions to a system of real polynomial equations.
This is of unquestionable importance in applications of mathematics.
Even the existence of real solutions is not guaranteed; oftentimes there are few
or no real solutions, and all complex solutions must be found to determine
if this is the case.
We give a method to construct families of polynomial systems that have
nontrivial lower bounds on their numbers of real solutions, guaranteeing
the existence of real solutions.

{\it Geometric problems} with lower bounds on their numbers of real solutions 
are a recent discovery.
Kharlamov and Degtyarev showed that of the 12 (a priori complex) rational cubics
passing through 8 real points in the plane, at least 8 are
real~\cite[Prop.~4.7.3]{DeKh00}.
This was generalized by Welschinger~\cite{We03},
Mikhalkin~\cite{Mi03,Mi04}, and Itenberg, Kharlamov, and Shustin~\cite{IKS03}, to 
rational curves passing through real points on toric surfaces.
Welschinger discovered an invariant which gives a lower bound, and work of
Mikhalkin and of Itenberg, Kharlamov, and Shustin shows that this
lower bound is non-zero and in fact quite large. 
If $N_d$ is the Kontsevich number of such complex rational
curves~\cite{KM94} and $W_d$ is Welschinger's invariant, then
$\log N_d$ and $\log W_d$ are each asymptotic to $3d\log d$.

At the same time, Eremenko and Gabrielov~\cite{EG01a,EG01b} computed the degree of the
Wronski map on the real Grassmannian of $k$-planes in $n$-space.
It is non-trivial when $n$ is odd.
This degree is a lower bound on the number of real solutions to certain problems
from the Schubert calculus on this Grassmannian.
In its formulation as a Wronski determinant, their work implies the existence of
many inequivalent $k$-tuples of polynomials of even degree having a given real
polynomial as their Wronskian. \smallskip

These results highlight the importance of developing a theoretical framework
to explain this phenomenon.
Our main purpose is to provide such a framework for sparse
polynomial equations.
We are inspired by the work of Eremenko and Gabrielov.
Our lower bound is the topological degree of a linear projection on
an oriented double cover of a toric variety.
In Section 1, we formulate a polynomial system as the fibers of a map from a
toric variety and define the characteristic of such a map to be the degree of the
map lifted to a canonical double cover.
This has the same equations, but is taken in
the sphere covering real projective space.
(One method used by Eremenko and Gabrielov was to lift the Wronski map to a
double cover of non-orientable Grassmannians.)
This characteristic is defined only if the smooth points of the double cover are
orientable.  
We give criteria for this to hold in Section 2. 
In Section 3, we show how to compute the degree for some maps 
by degenerating the double cover of the toric variety into a union of oriented
coordinate spheres and then determine the degree of the same projection on
this union of spheres.

This method does not work for all linear
projections of toric varieties. 
For toric varieties associated to the order polytope of a poset $P$, there
are natural Wronski projections with a computable characteristic when the poset
$P$ is ranked mod 2.
That is, the lengths of all maximal chains in $P$ have the same parity. 
In this case, the degree is the sign-imbalance of $P$---the difference between
the numbers of even and of odd linear extensions~\cite{Wh01,St02}.
This pleasing construction is the subject of Section 4.

Section 5 contains further examples of this theory.
Grassmannians admit flat sagbi degenerations to such toric
varieties~\cite[Ch.~11]{Sturmfels_GBCP}.
For these, the Wronski map coincides with a linear projection we study, and we
are able to recover the results of Eremenko and Gabrielov in this way.
This is the topic of Section 6.
\smallskip

In Section 7, we give alternative proofs of our lower bound for the order polytope
of a poset $P$, when $P$ is the incomparable union of chains of lengths
$a_1,\dotsc,a_d$.
We show that the Wronski polynomial system in this case is equivalent to finding
all factorizations $f(z)=f_1(z)\dotsb f_d(z)$, where $f(z)$ is a fixed
polynomial of degree $a_1+\dotsb+a_d$, and the factors $f_1(z),\dotsc,f_d(z)$
that we seek have respective degrees $a_1,\dotsc,a_d$.
This reformulation reveals the existence of a new phenomenon for real
polynomial systems.
Not only do each of these systems possess a lower bound on their number of
real solutions, but certain numbers of real solutions cannot occur.
That is, there are gaps in the possible numbers of real solutions to
these polynomial systems.

\section{Systems of Sparse Polynomials as Linear Projections}\label{S:System}

Let $F(t_1,t_2,\dotsc,t_n)$ be a real polynomial.
The exponent vector $m=(m_1,m_2,\dotsc,m_n)$ of a monomial
$t^m:=t_1^{m_1}t_2^{m_2}\dotsb t_n^{m_n}$ appearing in $F$ is a point
in the integer lattice $\Z^n\subset\R^n$.
The {\it Newton polytope} $\Delta\subset\mathbb{R}^n$ of a polynomial $F$ is the
convex hull of its exponent vectors.
We study real solutions to systems of real polynomial equations
 \begin{equation}\label{E:Sparse_System}
   F_1(t_1,\dotsc,t_n)\ =\
   F_2(t_1,\dotsc,t_n)\ =\   \dotsb\ =\
   F_n(t_1,\dotsc,t_n)\ =\ 0\,,
 \end{equation}
where the polynomials $F_i$ have real coefficients with the same Newton
polytope $\Delta$.
By Kushnirenko's Theorem~\cite{Ko75}, there are at most 
$V(\Delta):=n!\mbox{vol}(\Delta)$ solutions to~\eqref{E:Sparse_System}
in the complex torus $(\C^{\times})^n$ and this number is attained
for generic such systems.
We call this number $V(\Delta)$ the {\it normalized volume} of $\Delta$.
We shall always assume that our polynomial systems are generic in
that they have $V(\Delta)$ solutions in $(\C^{\times})^n$, each
necessarily of multiplicity one.

\begin{ex}\label{Ex:Hex1}
 Suppose that we have a system of two polynomial equations of the form
\[
  a_i + b_ix + c_iy + d_ixy +e_ix^2y + f_ixy^2 + g_ix^2y^2\ =\ 0,\qquad
  \textrm{for}\ i=1,2\,.
\]
 The monomials which appear correspond to the lattice points
 $(0,0)$, $(1,0)$, $(0,1)$, $(1,1)$, $(2,1)$, $(1,2)$, and $(2,2)$, whose convex
 hull is a hexagon.
\[
  \epsfxsize=2.4cm\epsffile{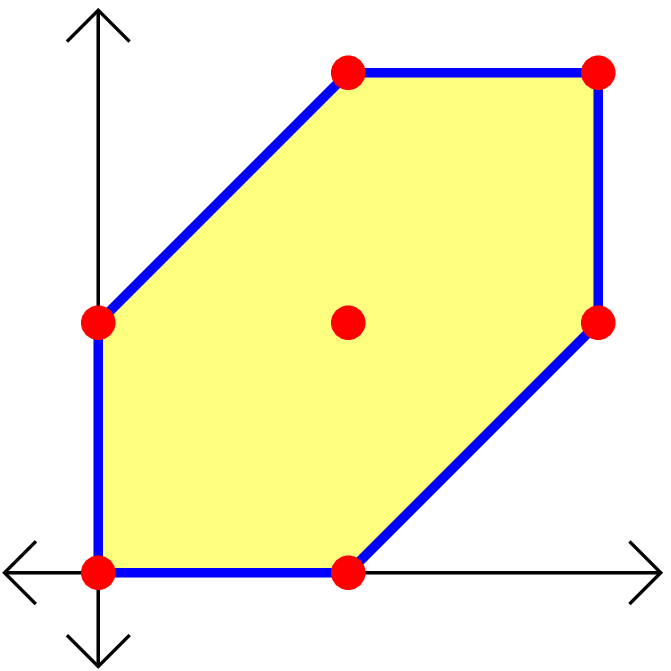}
\]
 This hexagon has Euclidean volume 3, and so we expect there to be
 $3\cdot 2! = 6$ complex solutions to this set of equations.
\end{ex}

The projective toric variety $X_\Delta$ associated to the polytope $\Delta$ is
the variety parametrized by the monomials in $\Delta$.
More precisely, let $\P^\Delta$ be the complex projective space
with coordinates $\{x_m\mid m\in\Delta\cap\Z^n\}$
indexed by the points of $\Delta\cap\Z^n$.
Then $X_\Delta$ is the closure of the image of the map
 \[
   \varphi_\Delta\ \colon\
    \begin{array}{rcl}
       (\C^\times)^n&\longrightarrow&\P^\Delta\\
       (t_1,t_2,\dotsc,t_n)&\longmapsto&
        [t^m\mid m\in\Delta\cap\Z^n]
    \end{array} \ .
 \]
This map is injective if and only if the affine span of $\Delta\cap\Z^n$ is
equal to $\Z^n$.

Linear forms on $\P^\Delta$ pull back along  $\varphi_\Delta$ to polynomials
with monomials from $\Delta\cap\Z^n$,
\[
   \varphi_\Delta^*\Bigr(\,
     \sum_{m\in\Delta\cap\Z^n} c_mx_m \Bigl)\ =\
   \sum_{m\in\Delta\cap\Z^n} c_mt^m\,.
\]
A system~\eqref{E:Sparse_System} of real polynomials with Newton polytope $\Delta$
corresponds to a system of $n=\dim X_\Delta$ real linear equations on $X_\Delta$,
that is, to the intersection of $X_\Delta$ with a real linear subspace $\Lambda$ of
codimension $n$ in $\P^\Delta$.

Let $E\subset\Lambda$ be a real hyperplane in $\Lambda$ disjoint from $X_\Delta$ -- a
linear subspace of $\P^\Delta$ complementary to $X_\Delta$.
Let $H (\simeq\P^n)\subset \P^\Delta$ be any real linear subspace of maximal
dimension $n$ disjoint from $E$.
Let $\pi_E$ be the linear projection with center $E$
 \begin{eqnarray*}
  \pi_E\ \colon\ \P^\Delta-E&\longrightarrow& H\\
    x&\longmapsto& \mbox{Span}(x,E)\cap H\,.
 \end{eqnarray*}
Then solutions to the system~\eqref{E:Sparse_System} correspond to
points in $X_\Delta\cap \pi_E^{-1}(p)$, where $p:=\pi_E(\Lambda)\in H$.

Set $Y_\Delta:=X_\Delta\cap \R\P^\Delta$, the real points of the toric variety
$X_\Delta$, and let $f$ be the restriction of $\pi_E$ to $Y_\Delta$.
We could also consider the closure of the image of $(\R^{\times})^n$
under $\varphi_\Delta$.
These objects coincide if and only if the restriction of
$\varphi_\Delta$ to $(\R^{\times})^n$ is injective, which occurs if
and only if the lattice spanned by $\Delta\cap\Z^n$ has odd index in
$\Z^n$.
We shall always assume that this index is odd.

Then real solutions to the system~\eqref{E:Sparse_System} are the elements
in the fiber $f^{-1}(p)$ of the linear projection $f$
\[
  f\ \colon\ Y_\Delta\ \subset\ \R\P^\Delta\
    \stackrel{\pi_E}{---\to}\ H_\R \simeq\R\P^n\,.
\]
If both $Y_\Delta$ and $\R\P^n$ are oriented, then the absolute value of the
topological degree of the map $f$ is a lower bound for the number of
points in $f^{-1}(p)$.
Our assumption on the genericity of the original system~\eqref{E:Sparse_System}
implies that $p$ is a regular value of the map $f$.

In general $Y_{\Delta}$ and $\R\P^n$ are not necessarily orientable.
Given a normal projective variety $Y\subset \R\P^N$ of dimension $n$, let
$Y^+\subset S^N$ be the subvariety of the sphere given by the same homogeneous
equations as $Y$.
Then $Y^+\to Y$ is a double cover.
Likewise, if $f\colon Y\to \R\P^n$ is the restriction of a linear projection
$\pi\colon \R\P^N-\to\R\P^n$ to $Y$, then we let
$f^+\colon Y^+\to S^n$ be the restriction of that projection lifted to
the corresponding spheres.
We obtain the commutative diagram, where the vertical arrows are 2 to 1 covering
maps.
\[
   \begin{picture}(160,50)(-5,0)
   \put( -3, 0){$f$}
   \put( 11, 0){$\colon$}
   \put( 23, 0){$Y$}
   \put( 45, 0){$\subset$}
   \put( 57, 0){$\R\P^N$}
   \put( 84, 0){$\stackrel{\pi_E}{---\to}$}
   \put(135, 0){$\R\P^n$}

   \put( -3,34){$f^+$}
   \put( 11,34){$\colon$}
   \put( 23,34){$Y^+$}
   \put( 45,34){$\subset$}
   \put( 62,34){$S^N$}
   \put( 84,34){$\stackrel{\pi_E^+}{---\to}$}
   \put(140,34){$S^n$}

   \put( 28,28){\vector(0,-1){15}}
   \put( 68,28){\vector(0,-1){15}}
   \put(144,28){\vector(0,-1){15}}
    \end{picture}
\]

\begin{defn}\label{D:char}
 Suppose that the manifold $Y^+_{\rm sm}$ formed by smooth points of $Y^+$ is
 orientable.
 Fix an orientation of $Y^+_{\rm sm}$ and define the {\it characteristic} of $f$,
 $\ch(f)$, to be the absolute value of the topological degree of
 $f^+\colon Y^+\to\R\P^n$.
 This does not depend upon the choice of orientation of $Y^+_{\rm sm}$ if
 it is connected.
 If $Y^+_{\rm sm}$ is not connected, then $\ch(f)$ could depend upon the choice
 of orientation of its different components.
 Since $Y$  is normal, the set of singularities $Y^{+}_{\rm sing}$ has codimension
 at least 2. Hence $\R\P^n\setminus\pi(Y_{\rm sing}^+)$ is connected
 and this notion is well-defined.
 \end{defn}

 Suppose that $Y$ is orientable. Consider the orientation on $Y^+$
 that is pulled back from 
 $Y$ along the covering map $S^N\to\RP^N$. 
If $\RP^n$ is not orientable then  $\ch(f)=0$.  
 If $\RP^n$ is orientable then the characteristic $\ch(f)$ is equal
 to the topological degree of $f$.

We record the obvious, fundamental, and important property of this notion. 

\begin{prop}\label{P:deg_leq_solutions}
 If $p\in\R\P^n$ is a regular value of $f$, then the number of points in a fiber
 $f^{-1}(p)$ is bounded below by its characteristic $\ch(f)$.
\end{prop} 

 According to Eremenko and Gabrielov~\cite{EG01b}, this notion is due to  
 Kronecker~\cite{Kr1968}, who defined the characteristic of a regular map
 $\R\P^2\to\R\P^2$ in this manner.
 Note that if $Y^+_{\rm sm}$ is not connected, then different choices of
 orientation of the components of $Y^+_{\rm sm}$ may give different values for
 $\ch(f)$.
 Each value for $\ch(f)$ is a lower bound on the number of points in a
 fiber $f^{-1}(p)$ above a regular value $p$ of $f$.
 Optimizing these choices is beyond the scope of this paper.

\section{Orientability of Real Toric Varieties}\label{S:orientibility}

The elementary definition of $Y_\Delta$ given in Section~\ref{S:System}, as the
real points of the variety parametrized by monomials in $\Delta\cap\Z^n$, is
inadequate to address the orientability of $Y_\Delta^+$. 
More useful to us is Cox's construction of $X_\Delta$ as a quotient
of a torus acting on affine space, as detailed in~\cite[Theorem 2.1]{Co95}.
Let $\Delta\subset\R^n$ be a polytope with vertices in the integer lattice
$\Z^n$ and suppose that it is given by its facet inequalities
 \begin{eqnarray*}
  \Delta&=&\{x\in\R^n\mid \calA\cdot x \geq -b\}\\
        &=&\{x\in\R^n\mid a_i\cdot x\geq -b_i, \ i=1,\dotsc,r\}\,,
 \end{eqnarray*}
where $a_i\in\Z^n$ is the primitive inward-pointing normal to the $i$th facet of
$\Delta$. 

\begin{ex}\label{Ex:Hex2}
 If $\Delta$ is the hexagon of Example~\ref{Ex:Hex1}, then
\[
   \calA\ =\  
     \left(\begin{array}{rrrrrr}0&-1&-1&0&1&1\\1&1&0&-1&-1&0\end{array}\right)^T
   \quad\mbox{ and }\quad
   b\ =\  \left(\begin{array}{rrrrrr}0&1&2&2&1&0\end{array}\right)^T\,.
\]
\end{ex}

Let $z=(z_1,\dotsc,z_r)\in\C^r$.
For each $m\in\Delta\cap\Z^n$, set
\[
   z(m)\ :=\ \prod_{i=1}^r z_i^{a_i\cdot m + b_i}\,,
\]
and consider the map $\psi_\Delta$ defined by
\[
  \psi_\Delta(z)\ =\ [z(m)\mid m\in\Delta\cap\Z^n]\ \in \ \P^\Delta\,.
\]
This map is undefined on the zero locus $B_\Delta$ of the monomial ideal
\[
  \langle z(m)\mid m\in\Delta\cap\Z^n\rangle.
\]
Note that $z_i$ appears in $z(m)$ if and only if $m$ does not lie on the $i$th facet.
Define the vertex monomial $z^{(v)}$ to be the product of all $z_i$ such that $v$ misses the $i$th facet.
Then $B_\Delta$ is the zero locus of the monomial ideal 
\[
  \langle z^{(v)}\mid v\ {\rm a\ vertex\ of}\ \Delta\rangle.
\]

The monomial $z^b=z_1^{b_1}\dotsb z_r^{b_r}$ divides each 
component $z(m)$ of $\psi_\Delta(z)$.
Removing these common factors from $\psi_\Delta(z)\in\P^\Delta$ 
shows that $\psi_\Delta$ factors through $\varphi_\Delta$, at least for
$z$ in the torus $(\C^\times)^r$.
For $z\in(\C^\times)^r$, we have $\psi_\Delta(z)=\varphi_\Delta\circ\phi_\Delta(z)$,
where 
\[
  \phi_\Delta\ \colon \ 
  (z_1,z_2,\dotsc,z_r)\ \longmapsto\ 
  (\dotsc, z_1^{a_{1i}}z_2^{a_{2i}}\dotsb z_r^{a_{ri}},\dotsc)\,.
\]
Since $\Delta$ has full dimension, this map is surjective and so the
image of $\psi_\Delta$ is dense in the projective toric variety $X_\Delta$.
Since $\varphi_\Delta$ is injective, two points of $(\C^\times)^r$ have the same
image under $\psi_\Delta$ if and only if they are equal modulo the kernel
$G_\Delta$ of $\phi_\Delta$ 
 \[
   G_\Delta\ :=\
   \{\mu\in(\C^\times)^r \mid  1=\prod_{i=1}^r \mu_i^{a_{ij}}\ 
      \mbox{for each }    j=1,\dotsc,n\}\,.
 \]

The map $\psi_\Delta$ almost identifies $X_\Delta$ as the
quotient of $\C^r-B_\Delta$ by $G_\Delta$.
The difficulty is that $G_\Delta$-orbits on $\C^r-B_\Delta$ are not necessarily
closed and so the geometric quotient $(\C^r-B_\Delta)/G_\Delta$ may not be 
Hausdorff. 
If $\Delta$ is a simple polytope (each vertex lies on exactly $n$
facets), then this does not occur and $X_\Delta$ is the geometric
quotient. 
In general, $X_\Delta$ is the closest variety to the non-Hausdorff quotient.
More precisely, it is the quotient in the category of schemes, the 
{\it categorical quotient}, written $(\C^r-B_\Delta)/\!/G_\Delta$.

\begin{prop}[Theorem~2.1~\cite{Co95}]
  Suppose that $\Delta\cap\Z^n$ affinely spans $\Z^n$.
  Then the abstract toric variety $X_\Sigma$ defined by the normal fan $\Sigma$
  of $\Delta$ is the categorical quotient $(\C^r-B_\Delta)/\!/G_\Delta$, and the
  map $\psi_\Delta$ induces an isomorphism of toric varieties
  $X_\Sigma \to X_\Delta$.
  This categorical quotient is a geometric quotient if and only if $\Delta$ is
  simple.
\end{prop}

If we restrict the map $\psi_\Delta$ to $\R^r-B_\Delta$, then
its image lies in the real toric variety $Y_\Delta$, but this image is not in
general equal to $Y_\Delta$. 

\begin{prop}
  The image of\/ $\R^r-B_\Delta$ under the map $\psi_\Delta$ is
  equal to $Y_\Delta$ if and only if the index of the lattice $\Lambda_\calA$
  spanned by the columns of $\calA$ in its saturation
  $\Lambda_\calA\otimes_\Z\Q$ is odd.
 \end{prop}

The lattice $\Lambda_\calA$ for the hexagon of Example~\ref{Ex:Hex1} is
saturated as $\calA$ has a $2\times 2$ minor with absolute value 1.

\noindent{\it Proof.}
 It suffices to show that the image of $(\R^\times)^r$ under the map
 $\phi_\Delta$ is equal to the real points $(\R^\times)^n$ of 
 $(\C^\times)^n$ if and only if  the index of $\Lambda_\calA$ in its saturation
 is odd.
 
 Invertible integer row and column operations reduce $\calA$ to its Smith normal
 form 
\[
  \left[
   \begin{matrix}
     a_1   &0     &\dots &0     \\
     0     &a_2   &\dots &0     \\
     \vdots&      &\ddots&      \\
     0     &0     &\dots &a_n   \\
     0     &0     &\dots &0     \\   
     \vdots&\vdots&\dots &\vdots\\
     0     &0     &\dots &0     
   \end{matrix}
  \right]\ .
\]
 These operations do not change the index of the lattice $\Lambda_\calA$ in
 its saturation. It follows that the index is equal to the product $a_1\cdots a_n$.
 The image of $(\R^\times)^r$ under the map
 $\phi_\Delta$ is equal to $(\R^\times)^n$ if and only if the map from
 $(\R^\times)^r$ to $(\R^\times)^n$ defined by 
 $(z_1,\dots,z_r)\mapsto (z_1^{a_1},\dots,z_n^{a_n})$ is surjective,
 which happens if and only if the product $a_1\cdots a_n$ is odd.
\hfill \raisebox{-3pt}{\epsfxsize=12pt\epsffile{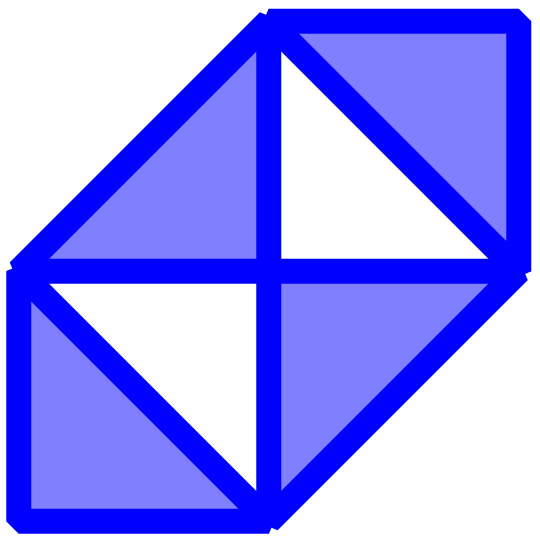}}\medskip

We address the orientability of $Y^+_\Delta$.
First, set $\ell:=\#\Delta\cap\Z^n$ and consider the map 
$g\colon \C^r\to\C^\ell$ which lifts the map
$\psi_\Delta\colon\C^r\to\P^{\ell-1}=\P^\Delta$
\[
  g\ \colon\ z\ \longmapsto\ (z(m)\mid m\in \Delta\cap\Z^n)\,.
\]
If we let $\gamma$ be the map from $G_\Delta$ to $\C^{\times}$ defined by
\[
  \gamma(\mu_1,\mu_2,\dotsc,\mu_r)\ =\ 
  \mu_1^{b_1}\mu_2^{b_2}\dotsb\mu_r^{b_r}\ =:\ \mu^b\,,
\]
then points $w,z\in(\C^{\times})^r$ have the same image in $\C^\ell$ if and only if
$wz^{-1}\in\mbox{ker}(\gamma)$. 
More generally, the fibers of the map
$\C^r-B_\Delta\to \C^\ell$ are unions of orbits of $\mbox{ker}(\gamma)$.

Let $\R_>$ be the positive real numbers.
The real cone over $Y_\Delta$ is $Z_\Delta:=(\R^\ell-\{0\})\cap g(\C^r)$.
Then the double cover $Y^+_\Delta$ of $Y_\Delta$ is the quotient 
$Z_\Delta/\R_>$. 
If we assume that the column space of $\calA$ has odd index in its saturation,
then there are two cases to consider.
 \begin{enumerate}
  \item[(i)] $Z_\Delta=g(\R^r)$, \, or 
  \item[(ii)] $Z_\Delta$ is the disjoint union of $g(\R^r)$ and $-1\cdot g(\R^r)$.
 \end{enumerate}
These cases are distinguished by the image of the map $\gamma$, when restricted
to the real points $G_\Delta(\R)$ of $G_\Delta$.
This also serves to describe $Y^+_\Delta$.
Set $K:=\gamma^{-1}(\R_>)$.

\begin{prop}
 With the above definitions, we have
 \begin{enumerate}
  \item[(i)] If $K\subsetneq G_\Delta(\R)$ so that 
             $\gamma(G_\Delta(\R))=\R^{\times}$, then
             $Z_\Delta=g(\R^r-B_\Delta)$, and so $Y^+_\Delta$ is the image of 
             $\R^r-B_\Delta$ under the composition
   \begin{equation}\label{E:composition}
       \R^r-B_\Delta\ \xrightarrow{\ g\ }\ 
       Z_\Delta\ \longrightarrow\ 
       Z_\Delta/\R_>\ =\ Y^+_\Delta\,.
   \end{equation}

  \item[(ii)] If $K=G_\Delta(\R)$, so that $\gamma(G_\Delta(\R))=\R_>$, then
              $Z_\Delta\neq g(\R^r-B_\Delta)$ but we have
     \[
           Z_\Delta\ =\  g(\R^r-B_\Delta)\coprod -g(\R^r-B_\Delta)\,.
     \]
              Furthermore,
              $Y^+_\Delta$ has two components, each isomorphic to
              $Y_\Delta$, and these components are interchanged by the antipodal
              map on the sphere $S^\Delta=S^{\ell-1}$, and one component is the
              image of\/ $\R^r-B_\Delta$ under the map~\eqref{E:composition}.
 \end{enumerate}
\end{prop}

We state our main result on the orientability of $Y^+_\Delta$.

\begin{thm}\label{T:orientability}
  Suppose that the lattice affinely spanned by $\Delta\cap\Z^n$ has odd index in
  $\Z^n$ and that $\Lambda_\calA$ has odd index in its saturation.
  If there is a vector $v$ in the integer column span of $[\calA:b]$, all of
  whose components are odd, then the standard orientation of $\R^r$ induces an
  orientation on the smooth part of $Y^+_\Delta$ via the map
  $\psi_\Delta\colon \R^r-B_\Delta\to Y^+_\Delta$.
\end{thm}

\begin{rem}\label{R:orient_Y} 
 If there is a vector $v$ in the integer column span of $\calA$, all of
  whose components are odd, then the orientation of $\R^r$
  induces an orientation on the smooth part of $Y_{\Delta}$. The proof
  of this statement is analogous to the proof of Theorem~\ref{T:orientability}.
\end{rem}

\begin{rem}\label{R:char_zero}
 In general, we may not know if either $Y_{\Delta}$ or $Y_{\Delta}^{+}$ are
 orientable. 
 The positive part $Y_{\Delta}^>$ of $Y_\Delta$ is the
 intersection of $Y_\Delta$ with the positive orthant of $\P^\Delta$ is always
 orientable, as it is isomorphic to $\Delta$, as a manifold with
 corners~\cite[\S 4]{Fu93}.
\end{rem}

When the hypotheses of Theorem~\ref{T:orientability} are satisfied, we 
assume that the smooth points of $Y^+_\Delta$ have the orientation
induced by $\psi_\Delta$, and we say that $Y^+_\Delta$ is 
{\it Cox-oriented}.
If $\Delta$ is the hexagon of Example~\ref{Ex:Hex1}, then $Y^+_\Delta$ is
Cox-oriented as 
it is smooth and the vector with all components 1 is the sum of the three
columns of the $6\times 3$-matrix $[\calA:b]$.\medskip

\noindent{\it Proof of Theorem~$\ref{T:orientability}$.}
 Recall that the subgroup $K\subset G_{\Delta}(\R)$ is
\[
  K\ :=\ \gamma^{-1}(\R_>)\ =\ 
   \{\mu\in G_{\Delta}(\R)\mid \mu^b>0\}\,.
\]
 We claim that if $\mu\in K$, then $\det(\mu)=\mu_1\mu_2\dotsb\mu_r>0$, so that
 $K$ preserves the standard orientation on $\R^r$.
 Indeed, let $c=(c_1,\dotsc,c_r)$ be an integer vector with each component
 $c_i$ odd such that $c-kb\in\Lambda_{\calA}$ for some $k\in\Z$.  
 Then  $\mu^c=(\mu^b)^k>0$, and so we have $\det\mu>0$, as each component of $c$ is
 odd (for then $\mu^c/\det(\mu)$ is a square). 

 Thus if $U\subset\R^r$ is an open subset with $K\cdot U=U$ such that 
 every orbit of $K$ is closed in $U$, then the smooth part of the
 quotient $U/K$ has an orientation induced by the standard
 orientation of $\R^r$.

 For each face $F$ of the polytope $\Delta$, let $\tau_F\subset\R^n$
 be the cone generated by the primitive inward-pointing normal
 vectors to the facets containing $F$---these generators are the
 rows of $\calA$ corresponding to the facets containing $F$.
 Set $U_F\subset\C^r$ to be the complement of the variety defined by
 the monomial ideal $\langle z(m)\mid m\in F\rangle$.
 This is the set of points $(z_1,\dotsc,z_r)\in\C^r$ such that
 $z_i\neq 0$ if the $i$th facet of $\Delta$ does not contain $F$.
 We have $G_\Delta\cdot U_F=U_F$.

 If the cone $\tau_F$ is simplicial, then every $G_\Delta$-orbit of
 $U_F$ is closed.
 The arguments that show this in the proof of Theorem 2.1
 of~\cite{Co95} show that the same is true of the $K$-orbits
 of $U_F(\R)$.
 Furthermore, $U_F/G_\Delta$ and also $U_F(\R)/K$ is smooth if the generators of
 $\tau_F$ in addition generate a saturated sublattice of $\Z^n$.

 If $F$ is a facet, then $\tau_F$ is just a ray generated by a primitive vector,
 and is thus simplicial.
 If we let $F$ run over the facets of $\Delta$, the quotients $U_F(\R)/K$ are
 glued together along the common subtorus, which is $U_\Delta(\R)$.
 As each piece and the torus is oriented by the canonical orientation of $\R^r$
 under the quotient by $K$, this union $W$ is a smooth and oriented subset of
 $Y^+_\Delta$.
 Moreover, the difference $\overline{W}-W$ has codimension 2 (this is part of
 the proof that toric varieties are smooth in codimension 1).

 Thus the image of $\R^r-B_\Delta$ in $Y^+_\Delta$ is smooth and oriented in
 codimension 1.
 This either is dense in $Y^+_\Delta$ or in one of the two isomorphic
 components of $Y^+_\Delta$.
 This completes our proof of the theorem.
\hfill \raisebox{-3pt}{\epsfxsize=12pt\epsffile{figures/QED.eps}}\medskip

\section{Computation of the characteristic}\label{S:char}
Let $\Delta\subset\R^n$ be a lattice polytope
with a regular triangulation $\Delta_{\omega}$ defined by
a lifting function $\omega:\Delta\cap\Z^n\to\Z_{\geq 0}$, and
$Y_{\Delta}\subset\P^{\Delta}$ the real toric variety parametrized by the
monomials in $\Delta$.
Assume that $\omega$ is convex, which means that all integer points of $\Delta$
are vertices of $\Delta_\omega$.
Call simplices with odd normalized volume {\it odd} and simplices with even
normalized volume {\it even}.

\begin{defn}
 A triangulation $\Delta_{\omega}$ is {\it balanced\/} if its vertex-edge
 graph is $(n+1)$-colorable.
 This means that there exists a map $\kappa$ from the integer points of
 $\Delta\cap\Z^n$ to the vertices of the standard simplex which is a bijection
 on each simplex in the triangulation $\Delta_{\omega}$.
 We call this map $\kappa$ a {\it folding} of $\Delta_{\omega}$.
\end{defn}

A triangulation is balanced if and only if its dual graph is bipartite.
For the direct implication, note that an orientation of the standard simplex
induces orientations of the simplices in the triangulation $\Delta_w$ via the
map $\kappa$.
This induced orientation changes when passing to an adjacent simplex.
The other implication is~\cite[Corollary~11]{J02}.

\begin{defn}
 For a balanced triangulation $\Delta_{\omega}$, assign $+$ or $-$ to each of
 the simplices so that every two adjacent simplices have opposite signs.
 Disregard even simplices and define the {\it signature}
 $\sigma(\Delta_{\omega})$ of $\Delta_{\omega}$ to be the absolute value of the
 difference of the numbers of odd simplices with $+$ and odd simplices with
 $-$.
\end{defn}

For each $m\in\Delta\cap\Z^n$, fix a nonzero real number $\alpha_m$ whose sign depends
only upon $\kappa(m)$.
Call this vector $(\alpha_m\mid m\in\Delta\cap\Z^n)$ a {\it weight function} for
$\Delta$.
As the vertices of the standard simplex are the standard basis
vectors in $\P^n$ and the vertices of the triangulation are the basis vectors of
$\P^{\Delta}$, the folding $\kappa$ defines a linear projection, called the
{\it Wronski projection} $\pi_\alpha\colon\P^{\Delta}\to\P^n$ sending each basis
vector $e_m$ of $\P^{\Delta}$ to $\alpha_m e'_{\kappa(m)}$, where $e'_i$ is a
basis vector of $\P^n$.
If $\alpha$ is constant, then we omit it from our system of notation as it has
no effect.

A linear form $\Lambda$ on $\P^n$ pulls back along
$\pi_\alpha$ to a polynomial $F$ of the from
\[
   F\ =\ \sum_m c_{\kappa(m)}\, \alpha_m x^m\ ,
\]
where $c_0,c_1,\dotsc,c_n\in\R$.
We call such a polynomial a Wronski polynomial for the triangulation
$\Delta_w$ and the weight function $\alpha$.
A {\it Wronski polynomial system} for the triangulation $\Delta_w$ and
weight function $\alpha$ is a system of $n$ such polynomials, all with
weight function $\alpha$.
Solutions to such a Wronski system correspond to a fiber of the Wronski map
$\pi_\alpha$.

\begin{ex}\label{Ex:Hex_fold}
The lattice hexagon of Example~\ref{Ex:Hex1} admits a regular unimodular balanced
triangulation induced by the lifting function taking value 0 at its center
$(1,1)$, 3 at the vertices $(0,0)$ and $(2,2)$, and 1 at the remaining 4
vertices.
This triangulation defines a 3-coloring of the vertices indicated
by the labels $a$, $b$, $c$ in the Figure~\ref{F:HEX_FOLD}.
 \begin{figure}[htb]
  \[
   \begin{picture}(120,100)(0,5)
     \put(10,9){\epsfxsize=1.2in\epsffile{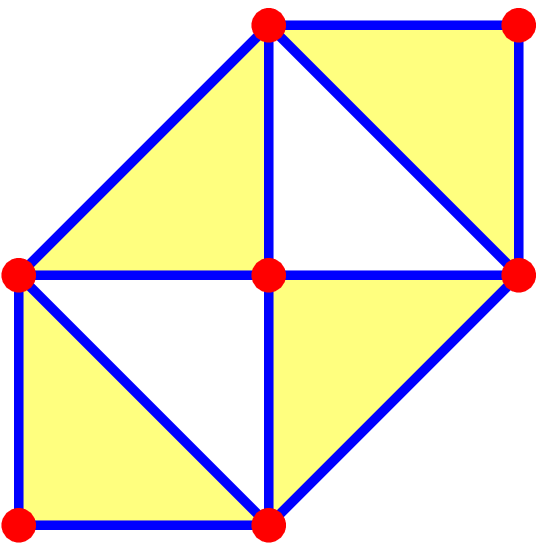}}
                      \put(46,97){$b$} \put(95, 97){$a$}
     \put(2,51){$c$}  \put(57, 56){$a$}   \put(98, 51){$c$}
     \put(4, 1){$a$}  \put(50,  0){$b$}
    \end{picture}
     \qquad \qquad
    \begin{picture}(50,50)(0,0)
    \put(10,9){\epsfxsize=1.2in\epsffile{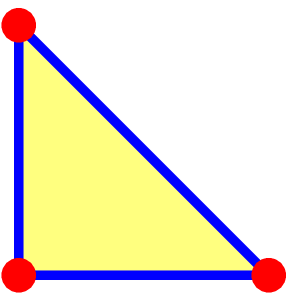}}
    \put(2,51){$c$}
    \put(4, 1){$a$}  \put(50,0){$b$}
    \end{picture}
  \]
 \caption{}
 \label{F:HEX_FOLD}
 \end{figure}
We illustrate the bipartite dual graph by shading the positive simplices.
For the constant weight function, this defines a Wronski projection from $\P^6$
to $\P^2$ by
\[
   [x_{00},x_{10},x_{01},x_{11},x_{21},x_{12},x_{22}]\ \longmapsto\
   [x_{00}+x_{11}+x_{22}, x_{10}+x_{12}, x_{01}+x_{21}].
\]
The corresponding Wronski polynomials have the form
 \[
   a(1+xy+x^2y^2) + b(x+xy^2) + c(y+x^2y)\ =\ 0\,,
 \]
 where $a$, $b$, and $c$ are arbitrary real numbers.
\end{ex}

The lifting function $\omega=(\omega_m\colon m\in\Delta\cap\Z^n)$ defines a
partial term order on the coordinate ring $\R[x_m\mid m\in\Delta\cap\Z^n]$
of $\P^{\Delta}$
by $x^b>x^c$ if $b\cdot\omega<c\cdot\omega$, where $b\cdot\omega$ and
$c\cdot\omega$ denote the standard scalar product.
It also defines an action of $\R^{\times}$ on $\P^{\Delta}$ by
 \begin{equation}
  s. x_m\ =\ s^{-\omega_m}\cdot x_m\,.
 \end{equation}
The corresponding action on $\R[x_m\mid m\in\Delta\cap\Z^n]$ is the dual action
$$s. g(x)=g(s^{-1}. x)$$
Thus a monomial $x^b$ is transformed into $s^{b\cdot\omega}x^b$.
The monomials in the initial form  $\ini_{\omega}g$ of $g$ are multiplied by
the same power of $s$ in $g(x)$, which is less than the power of $s$ for the
other monomials.
Dividing $s. g$ by this lowest power $s^{b\cdot\omega}$ of $s$ we see that
 \[
   \lim_{s\to 0}s^{b\cdot\omega}s. g(x)\ =\ \ini_{\omega}g(x)\,.
 \]

Consider this for $Y_{\Delta}$.
The ideal $I(s. Y_{\Delta})$ of $s. Y_{\Delta}$ is
\[
   I(s. Y_\Delta)\ =\ \left\{ s. g(x) \mid g\in I(Y_\Delta)\right\}\,.
\]
Let $\ini(Y_{\Delta})$ be the variety defined by the initial ideal
$\ini_\omega I(Y_\Delta)$.
These arguments show that it is the scheme-theoretic limit of the family
$s.Y_\Delta$,
\[
  \lim_{s\to 0}s.Y_\Delta\ =\ \ini(Y_{\Delta})\,.
\]
If we define $s.\alpha$ by $(s.\alpha)_m=s^{m\cdot\omega}\alpha_m$, then the
Wronski map $\pi_\alpha$ on $s.Y_\Delta$ is equivalent to the Wronski map
$\pi_{s.\alpha}$ on $Y_\Delta$.

The family $\{s.Y_{\Delta}\mid s\in(0,1]\}\cup\ini(Y_\Delta)$ is a
{\it toric degeneration} of $Y_{\Delta}$.
This action lifts to the sphere, giving the family
$\{s.Y^{+}_{\Delta}\mid s\in(0,1]\}\cup\ini(Y_\Delta^+)$
in which $s.Y^{+}_\Delta$ and $\ini(Y_\Delta^+)$ are the subvarieties of the
sphere $S^{\Delta}$ given by the same homogeneous
equations as $s.Y_{\Delta}$ and $\ini(Y_\Delta)$.

By Kushnirenko's Theorem~\cite{Ko75} the number of complex solutions of the
Wronski system is equal to the normalized volume
$V(\Delta)$, which is an upper bound for the number of real solutions.
When the triangulation $\Delta_w$ is unimodular, Sturmfels~\cite{St94b} used
these toric degenerations to show that this upper bound is attained.
We use a toric degeneration to compute $\ch(f)$ which is by
Proposition~\ref{P:deg_leq_solutions} a lower bound for the number of
real solutions of a Wronski polynomial system.

\begin{thm}\label{T:computedegree}
 Suppose that the toric degeneration of $Y_\Delta$ does not meet the center of
 the Wronski projection $\pi_\alpha$ and $Y_{\Delta}^{+}$ is Cox-oriented.
 Then $\ch(f)$ is equal to the signature $\sigma(\Delta_{\omega})$
 of the triangulation $\Delta_{\omega}$.
 Moreover, if $s_0\in\R_>$ is minimal such that $s_0.Y_\Delta$ meets the
 center of projection, then
 $\ch\bigl(\pi_\alpha|_{s.Y_\Delta}\bigr)=\sigma(\Delta_{\omega})$, for any
 $0<s<s_0$.
\end{thm}

\begin{ex}\label{Ex:Hexagon}
 We observed that if $\Delta$ is the hexagon Example~\ref{Ex:Hex1}, then
 $Y^+_\Delta$ is Cox-oriented.
 In Example~\ref{Ex:Hex_fold}, we saw that $\Delta$ has
 a regular unimodular balanced triangulation (illustrated in
 Figure~\ref{F:HEX}), which has a signature of 2.
 The Wronski polynomials for the family $s.Y_\Delta$ have the form
 \begin{equation}\label{E:Hexs}
    a(s^2+xy+s^2x^2y^2) + bs(x+xy^2) + cs(y+x^2y)\ =\ 0\,,
 \end{equation}
 where $a$, $b$, and $c$ are arbitrary real numbers.
 The coefficients (in $s,x,y$) of $a,b$, and $c$ vanish where  $s.Y_\Delta$
 meets the center of projection.
 There are no real values of $x,y$ where these coefficients vanish for $s\neq 0$.
 Thus no variety $s.Y_\Delta$ in the family induced by the weight
 function meets the center of projection.

 By Theorem~\ref{T:computedegree}, for any given $s\neq 0$, two general
 polynomial equations of the form~\eqref{E:Hexs} will have at least 2 common
 real solutions.
 Figure~\ref{F:HEX} shows the two curves given by equations of the
 form~\eqref{E:Hexs} when $s=1$ with coefficients $(a,b,c)$ equal to
 $\Blue{(3,5,1)}$ and to $\Red{(1,-2,-3)}$, which meet in the two points
 indicated.
\begin{figure}[htb]
 \[
   \begin{picture}(120,108)
     \put(10,9){\epsfxsize=1.2in\epsffile{figures/Hex-triangulation.eps}}
                      \put(43,100){$xy^2$} \put(98, 99){$x^2y^2$}
     \put(0,51){$y$}  \put(58, 58){$xy$}   \put(100, 51){$x^2y$}
     \put(3, 0){$1$}  \put(50,  0){$x$}
    \end{picture}
     \qquad \qquad
   \epsfxsize=1.5in\epsffile{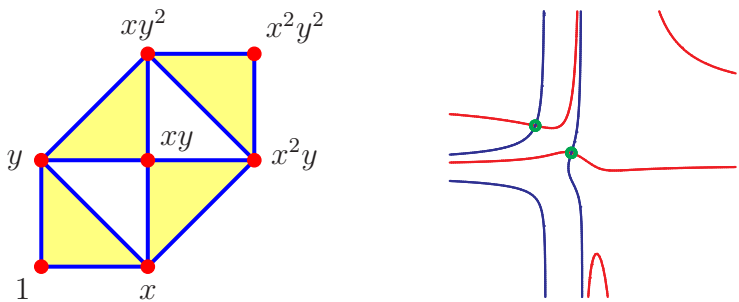}
 \]
 \caption{Hexagonal system}
 \label{F:HEX}
 \end{figure}
 We computed one million random instances of this polynomial system
 with $s=1$.
 Each one had exactly 2 real solutions.

 These computations, like all computations reported here, were done
 purely symbolically.
 The computation procedure involved generating random polynomial
 systems, and then computing a univariate eliminant for each system.
 This eliminant has the property that its number of real solutions
 equals the number of real solutions to the original system.
 This part of the computation was done with the computer algebra
 system Singular~\cite{Singular}.
 For all computations, except those reported in the last paragraph,
 the number of real roots for the eliminant were determined using
 an implementation of Sturm sequences in Singular.
 That implementation is inefficient for polynomials of degree 30, so
 the last computations in this paper used Maple's {\tt realroot}
 routine to compute the number of real solutions.
 
 We also computed 500,000 instances of the system~\eqref{E:Hexs} for $s\in (0,1)$.
 Of these, 429,916 had 2 real solutions 70,084 had 6 real solutions, and none had 4
 solutions.
 More precisely, $1000\cdot s$ was an integer chosen uniformly in $[1,999]$ and the
 coefficients $a,b,c$ were chosen uniformly in $[-60,60]$.

 Given fixed weights $\alpha_m$ for $m \in\Delta\cap\Z^2$, a Wronski
 polynomial with these weights is
\[
    a(\alpha_1+\alpha_{xy}xy+\alpha_{x^2y^2}x^2y^2)
  + b(\alpha_{x}x+\alpha_{xy^2}xy^2)
  + c(\alpha_{y}y+\alpha_{x^2y}x^2y)\ =\ 0\,.
\]
 We computed instances of such Wronski systems with 2, 4, or 6
 real solutions.
\end{ex}

\noindent{\it Proof of Theorem$~\ref{T:computedegree}$.}
We can assume that $\Delta_\omega$ has at least one odd simplex, for otherwise
the lower bound is trivial. Write $\pi$ for the Wronski projection $\pi_\alpha$.
It lifts to $\pi^+\colon S^{\Delta}\to S^n$ given by the same equations as
$\pi$.
Let $f_s^+$ be the restriction of $\pi^{+}$ to $s.Y_{\Delta^+}$ for $s\in(0,s_0)$
and $f_0^+$ the restriction of $\pi^+$ to $\ini(Y_{\Delta}^+)$.
Since the toric degeneration of $Y_{\Delta}$ does not meet the center of
projection $\pi$, the characteristic $\ch(f)$ is equal to the characteristic of
$f_s^{+}\colon s.Y^+_{\Delta}\to S^n$ for any $s\in (0,s_0)$.

\begin{figure}[htb]
  \[
   \begin{picture}(330,160)(0,5)
    \put(0,0){\includegraphics[width=4in]{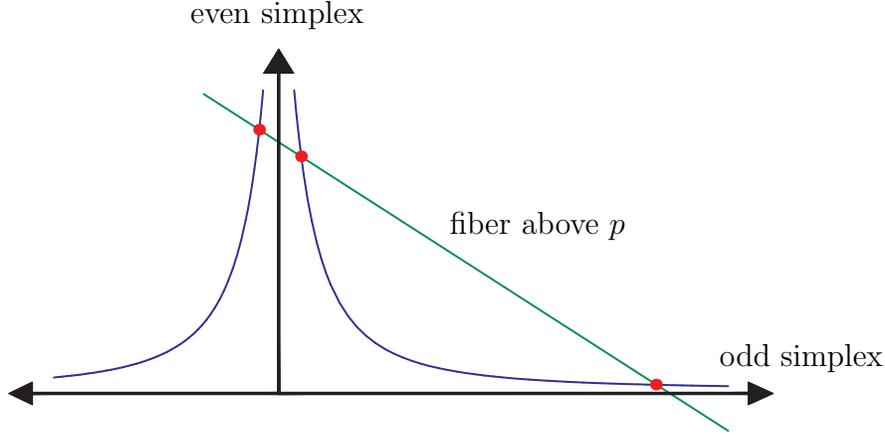}}
    \put(70,155){even simplex}
    \put(270,25){odd simplex}
    \put(168,75){fiber above $p$}
   \end{picture}
  \]
  \caption{Preimages near coordinate spheres.}\label{F:odd_even}
 \end{figure}

It is proved in Chapter 8 of \cite{Sturmfels_GBCP} that
\[
    {\rm Rad}(\ini_\omega I(Y_\Delta))\ =\
   \bigcap_{\tau}\langle x_m\,|\, m\not\in\tau\rangle\,,
\]
where the intersection is taken over all simplices $\tau$ of $\Delta_\omega$.
Thus $\ini_\omega(Y^+_{\Delta})$ is the union of coordinate $n$-planes
$\R\P^\tau$, one for each simplex $\tau$ in $\Delta_\omega$
 and $\ini_\omega(Y^+_{\Delta})$ is a similar union of coordinate $n$-spheres
$S^\tau$.
Thus a point $p\in S^n$ with non-zero coordinates
has one preimage $a_\tau$ under $f_0^+$ on each sphere
$S^\tau$.
The preimages of $p$ under
$f_s^+$ on $Y_s^+$ for small $s$ are clustered around
these $\{a_\tau\,\mid\,\tau\in\Delta_\omega\}$.
This is illustrated in Figure~\ref{F:odd_even}.
The \Red{preimages} are the \Red{dots}, the
\ForestGreen{linear subspace $(f^+)^{-1}(p)$} is the \ForestGreen{line}, and
the \Blue{toric variety $Y^+_\Delta$} is the \Blue{curve}.

When $s$ is small, consider the contribution to the characteristic of $f_s$ to
the solutions near $a_\tau$.
In a neighborhood of the point $a_\tau$ the projection $f_s^+$ is homotopic to
the coordinate projection $\pi_\tau$ to $S^\tau$ and therefore
we compute this local contribution using $\pi_\tau$.

This is easiest when $\tau$ is an even simplex, as in that case the restriction
$\pi_\tau|_{Y_\Delta^+}$ is not surjective and therefore this contribution to the
characteristic of $f^+_s$ is zero.
To see this, it is best to consider this projection in $\RP^\Delta$.
The composition
\[
  (\R^\times)^n\ \xrightarrow{\ \varphi_\Delta\ }\ \RP^\Delta\
  \xrightarrow{\ \pi_\tau\ }\ \R\P^\tau\ \simeq\ \RP^n\
\]
is the parametrization $\varphi_\tau$ of $\RP^\tau$ by the monomials
corresponding to integer points of $\tau$.
Since the affine span of the lattice points in $\Delta$ has odd index in $\Z^n$,
the map $\varphi_\Delta$ is an isomorphism between $(\R^\times)^n$ and the
dense torus in $Y_\Delta$.
Thus the restriction $\pi_\tau|_{Y_{\Delta}}$ is surjective if and
only if $\varphi_\tau$ maps $(\R^\times)^n$ onto the dense torus in
$\RP^n$. 
But this is not the case, as the integer points in $\tau$ span a sublattice of
$\Z^n$ with even index. For an odd simplex $\tau$, the map
$\varphi_\tau:(\R^{\times})^n\to(\R^{\times})^n$ is an isomorphism 
and therefore the degree of $\pi_\tau$ is 1. 

Pick a point $p=(p_0,\dots,p_n)$ in $S^n$ such that
$\sign p_i=-\sign \kappa(m)$ whenever $\kappa(m)=i$ 
where $\kappa$ is the folding of $\Delta_{\omega}$.
Then for each odd simplex there exists a unique preimage of $p$
 under $\pi_\tau$ and all of its components are positive. 
For an even simplex, there is an even number of  preimages with one
 of them having all components positive. 

Orient each of the coordinate spheres $S^\tau$ pulling back the
orientation of $S^n$ along $\pi$. 
Each of these orientations induces an orientation of the positive
part of  $Y_{\Delta}^+$. 
It remains to compare these induced orientations.

Consider two adjacent simplices in $\Delta_{\omega}$.
Let the vertices of the common facet be indexed by the variables $x_1,\dots,x_n$,
and the remaining two vertices by $x_0$ and $x_{n+1}$.
Then $x_0^{a_0}x_{n+1}^{a_{n+1}}=x_1^{a_1}\cdots x_n^{a_n}$ for some integers
$a_0,\dots,a_{n+1}$ with $a_0$ and $a_{n+1}$ positive.
Projections to the coordinate spheres of each simplex give local
coordinate charts for $Y^+_\Delta$, namely
$x_0,x_1,\dotsc,x_n$ and $x_{n+1},x_1,\dotsc,x_n$.
The Jacobian matric for this change of coordinates has the form
\[
   \left[\begin{matrix}
     \frac{\partial x_{n+1}}{\partial x_0}&0&\dotsb&0\\
     *&1&\dotsb&0\\
     \vdots&\vdots&\ddots&\vdots\\
    *&0&\dotsb&1\end{matrix}\right]\ ,
    \quad\mbox{{\rm where }}\ 
    \frac{\partial x_{n+1}}{\partial x_0}\ =\ 
    -\frac{a_0}{a_{n+1}}\frac{x_{n+1}}{x_0}\,.
\]
Since the Jacobian determinant is negative, these two charts belong
to different orienting atlases, and we need to count the
corresponding preimages with opposite signs. Therefore, $\ch(f)$ is
equal to 
the signature $\sigma(\Delta_\omega)$.
\hfill
\raisebox{-3pt}{\epsfxsize=12pt\epsffile{figures/QED.eps}}\medskip

\begin{rem}
Notice that this computation does not depend on the choice of
orientation of different connected components of the smooth part of
$Y_{\Delta}^+$. 
This implies that $\ch(f)=0$ whenever the smooth part of
$Y_{\Delta}^{+}$ is not connected.  
In particular, $\ch(f)=0$ if $Y_\Delta^+$ is isomorphic to two
copies of $Y_{\Delta}$. 
We have noted before that if $Y_{\Delta}$ is orientable but $\RP^n$
is not then $\ch(f)=0$ if the orientation on $Y_{\Delta}^+$ is
pulled back from $Y_{\Delta}$. 
We have proved that this last assumption is redundant: if
$Y_{\Delta}$ is orientable but $\RP^n$ is not then $\ch(f)=0$.
\end{rem}

\begin{lemma}\label{L:pos_comp}
If $\Delta_\omega$ contains only odd simplices and the sign of $\alpha_m$
depends only upon $\kappa(m)$, then there exists a regular value in
$z\in S^n$ all of whose preimages in $Y^+_\Delta$ under $f$ has all components
positive.
\end{lemma}

\noindent{\it Proof. }
Since $f(s^{-1}.x)=s^{-1}.z$ whenever $f_s(x)=z$, the statement follows from above.
\hfill \raisebox{-3pt}{\epsfxsize=12pt\epsffile{figures/QED.eps}}\medskip

 While in general we may not know if either $Y_{\Delta}$ or $Y_{\Delta}^{+}$ are
 orientable, the topological degree $\ch(f_>)$ of 
 $f_>:=f|_{Y_{\Delta}^>}\colon Y_{\Delta}^>\to f(Y_\Delta^>)$ is
 always well defined, as $Y^>_\Delta$ is orientible. 

\begin{cor}\label{C:Y_positive_char}
 Suppose that the toric degeneration of $Y_\Delta$ does not meet the center of
 projection $\pi$. Then $\ch(f_>)$ is equal to $\sigma(\Delta_{\omega})$.
\end{cor}

\section{Toric varieties from posets}\label{S:OrderPolytope}

Let $P$ be a finite partially ordered set (poset) with $n$ elements. 
We recall some definitions from the paper of Stanley~\cite{St86b}.

\begin{defn}
 {\it The order polytope} $O(P)$ of a finite  poset $P$ is the set of 
 points $y$ in the unit cube $[0,1]^P$ such that $y_a\leq y_b$ whenever  $a\leq
 b$ in $P$. 
\end{defn}

The vertices of the order polytope are the characteristic functions of (upper)
order ideals of $P$. 
Let $\calJ(P)$ be the set of such order ideals of $P$.
The {\it canonical triangulation} of the order polytope $O(P)$ is defined by the
linear extensions  (order-preserving bijections) of the poset $P$. 
Suppose that $P$ has $n$ elements and let $\lambda\colon P\to [n]$ be a linear
extension of $P$. 
For each $k=1,\dots,n$, let $a_k$ be the element of $P$ such that $\lambda(a_k)=k$.
Then $\lambda$ defines an
$n$-dimensional simplex $\tau_\lambda\subset O(P)$ consisting of  all $y$
satisfying 
\[
  0\leq y_{a_1}\leq\dots\leq y_{a_n}\leq 1\,
\]
The $\tau_\lambda$ are the simplices in a unimodular triangulation of $O(P)$. 
It is balanced as the association of an order ideal to its number of elements is
a proper coloring of its vertex-edge graph.
We will show in Lemma~\ref{L:canon_reg} that this triangulation is regular.

Fixing one linear extension of $P$ identifies each linear extension of $P$ with
a permutation of $P$, where the fixed extension is identified with the identity
permutation.
The sign of a linear extension is the sign of the corresponding permutation.

\begin{defn}
 The {\it sign imbalance} $\sigma(P)$ of a poset $P$ is
 the absolute value of the difference between the numbers of the
 positive and negative linear extensions of $P$.  
 If $\sigma(P)=0$ we say that $P$ is sign-balanced.
 Stanley studied this notion of sign-balanced posets~\cite{St02}.
\end{defn}

For an order ideal $J$, let $t^J:=\prod_{a\in J}t_a$ be
the monomial in $\R[t_a\mid a\in P]$ whose exponent vector is the vertex of $O(P)$
corresponding to the order ideal $J$. 
Let $|J|$ be the number of elements in the order ideal $J$. 
Fix a system of weights $\{\alpha_J\in\R^\times\mid J\in\calJ(P)\}$.
This gives the Wronski projection $\pi_\alpha$, Wronski polynomials,
and Wronski polynomial systems as in Section~\ref{S:char}.
Wronski polynomials for $\pi_\alpha$ have the form
 \[
      \sum_{J\in\calJ(P)} c_{|J|} \alpha_J t^J\,,
 \]
where  $c_0,\dotsc,c_{|P|}\in\R$.

\begin{thm}\label{T:order_poly_degree}
 Suppose that a finite poset $P$ is ranked mod $2$.
 For any choice $\alpha$ of weights, a Wronski polynomial system for the
 canonical triangulation of the order polytope of $P$ with weight $\alpha$ 
 has at least $\sigma(P)$ real solutions.
\end{thm}

The set $\calJ(P)$ of order ideals, ordered by inclusion, forms a ranked
distributive lattice $\calJ(P)$. 
Equations for the toric variety $Y_{O(P)}$ parametrized by the monomials in the
order polytope are nicely described by this lattice.
The lattice operations are $J\vee K=J\cap K$ and $J\wedge K=J\cup K$.
Its ideal $I$ is the Hibi ideal of this lattice~\cite{Hibi}
 \begin{equation}\label{E:toric_from_lattice}
   I\ =\ \langle x_J x_K\ -\  x_{J\wedge K}x_{J\vee K}\ \mid\ 
     J,K\in\calJ(P)\mbox{\rm \ are incomparable}\rangle\,.
 \end{equation}
The geometry of toric varieties associated to distributive lattices is discussed
in~\cite{Wa96}. 
Maximal chains of $\calJ(P)$ are the linear extensions of $P$. 
If two maximal chains differ by one element, they have opposite signs. 
Then $\sigma(P)$ is the absolute value of the difference between the number of
the positive maximal chains and the number of negative maximal chains.
We also call $\sigma(P)$ the sign-imbalance of the lattice $\calJ(P)$.

\begin{ex}\label{Ex:P2P2}
The toric variety $\P^2\times\P^2$ is defined by the order polytope of the poset
with the Hasse diagram
\[
  \begin{picture}(90,58)(-27,0)
    \put(-27,23){$P\ =$}
    \put(14,3){\epsffile{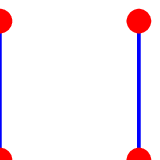}}
    \put(-2,0){$t_1$}   \put(61,0){$u_2$}
    \put(-2,49){$t_2$}  \put(61,49){$u_1$}
  \end{picture}
\]
Figure~\ref{F:P2P2} shows its lattice $\calJ(P)$ of order ideals and the six
maximal chains in $\calJ(P)$.
\begin{figure}[htb]
\[\setlength{\unitlength}{1.4pt}
  \begin{picture}(125,125)(-6,-1)
   \thicklines
               \put(33,112){$\{t_1,t_2,u_1,u_2\}$}
   \put(35,93){\Blue{\line(1,1){15}}}       \put( 78,93){\Blue{\line(-1,1){15}}}
   \put(8,84){$\{t_1,t_2,u_2\}$}          \put( 68,84){$\{t_2,u_1,u_2\}$}
   \put( 7,65){\Blue{\line(1,1){15}}}    \put(105,65){\Blue{\line(-1,1){15}}}
   \put(50,65){\Blue{\line(-1,1){15}}}   \put( 63,65){\Blue{\line(1,1){15}}}
   \put(-10,56){$\{t_1,t_2\}$}\put(41,56){$\{t_2,u_2\}$}\put(96,56){$\{u_1,u_2\}$}
   \put(22,37){\Blue{\line(-1,1){15}}}   \put( 90,37){\Blue{\line(1,1){15}}}
   \put(35,37){\Blue{\line(1,1){15}}}   \put( 78,37){\Blue{\line(-1,1){15}}}
   \put(21.5,28){$\{t_2\}$}              \put(75.5,28){$\{u_2\}$}
   \put(50, 9){\Blue{\line(-1, 1){15}}}   \put( 63, 9){\Blue{\line(1,1){15}}}
               \put(55, 0){$\emptyset$}
  \end{picture}
  \setlength{\unitlength}{1pt}
\]
\[
   \epsfxsize=4.in\epsffile{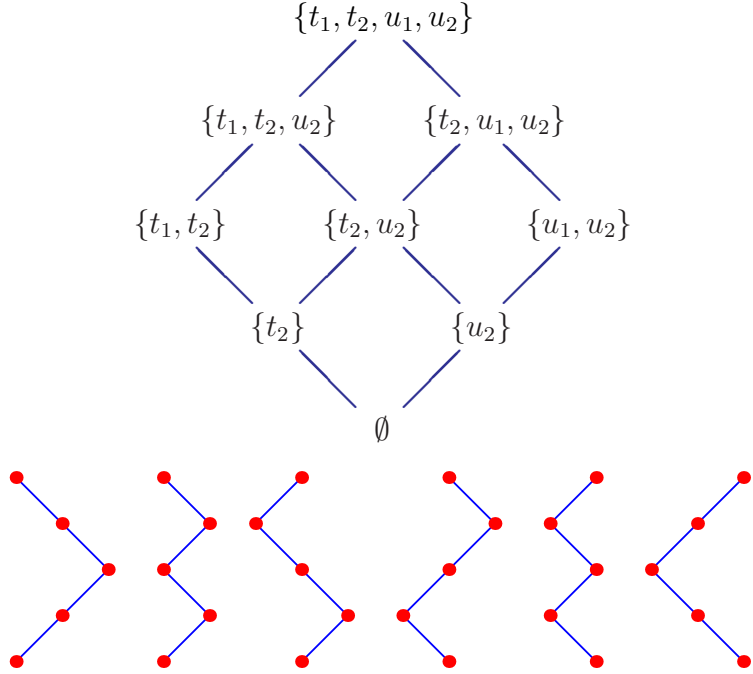}
\]
 \caption{The distributive lattice $\calJ(P)$ and its maximal chains.}
 \label{F:P2P2}
\end{figure}
The corresponding signs of these maximal chains are $+,-,+,+,-,+$, and
so the sign imbalance $\sigma(P)$ is 2. 
By Theorem \ref{T:order_poly_degree}
a generic system of 4 real equations of the form
\[
   c_0+c_1(t_2+u_2)+c_2(t_1t_2+t_2u_2+u_1u_2)+c_3(t_1t_2u_2+t_2u_1u_2)
   +c_4t_1t_2u_1u_2\ =\ 0\,,
\]
has at least 2 real solutions.
(For simplicity, the weights are constant, $\alpha_J=1$.)
\end{ex}

The following four lemmas reduce Theorem~\ref{T:order_poly_degree} to Theorem
\ref{T:computedegree}.

\begin{lemma}\label{L:canon_reg}
  The canonical triangulation of the order polytope $O(P)$ is regular.
\end{lemma}

\noindent{\it Proof.} 
 Define a lifting function $w(J)$ for each order ideal $J\in\calJ(P)$ by
\[
   \omega(J)\ :=\ -|J|^2\,.
\]
 For each simplex $\tau$ in the canonical triangulation we give a linear function
 $\Lambda$ on $\R^P$ such that $\Lambda(m)+\omega(m)\geq 0$ for all vertices of
 $O(P)$, with equality if and only if $m\in\tau$.

 For a linear extension $\lambda$, the vertices
 of the simplex $\tau_\lambda$ are
\[
  m_0\, :=\, (0,\dots,0),\quad\mbox{\rm and}\quad 
  m_k\, :=\, \sum_{i=1}^k e_{\lambda^{-1}(n-i+1)},
  \quad\mbox{\rm for}\quad  k=1,\dots,n\,, 
\]
 where $\{e_a\mid a\in P\}$ is the standard basis for $\R^P$. 
 The value of the lifting function at the vertex $m_k$ is $-k^2$.
 Define the linear function $\Lambda$ on $\R^P$ by
 $\Lambda(e_{\lambda^{-1}(n-k+1)})=2k{-}1$.
 Then $\Lambda$ is unique linear function on $\R^n$ such that 
 $\Lambda(m_k)+\omega(m_k)=0$.
 Indeed, 
\[
  \Lambda(m_k)\ =\ \sum_{i=1}^k \Lambda(e_{\lambda^{-1}(n-i+1)})
              \ =\ \sum_{i=1}^k 2i-1\ \ =\ k^2\ =\ -\omega(m_k)\,.
\]

 This also shows that if a vertex $m$ of $O(P)$ corresponds to an order ideal
 with $k$ elements, then $\Lambda(m)\geq k^2$ with equality only when $m=m_k$.
 Thus a vertex $m$ of $O(P)$ does not lie in $\tau_\lambda$ exactly when 
 $\Lambda(m)+\omega(m)>0$.
\hfill \raisebox{-3pt}{\epsfxsize=12pt\epsffile{figures/QED.eps}}\medskip

\begin{lemma}\label{L:Order_balanced}
  The canonical triangulation of the order polytope $O(P)$ is balanced. 
  Its signature is $\sigma(P)$, the sign imbalance of $P$.
\end{lemma}

\noindent{\it Proof.} 
 The folding map $m\mapsto |m|$ shows that the canonical
 triangulation is balanced.
 Linear extensions corresponding to adjacent simplices differ by a
 transposition and thus have opposite signs.
 The second statement is immediate.
\hfill \raisebox{-3pt}{\epsfxsize=12pt\epsffile{figures/QED.eps}}\medskip

\begin{lemma}\label{L:center_proj_orderpoly}
  For any choice of weights $\{\alpha_J\mid J\in\calJ(P)\}$, the toric
  degeneration of\/ $Y_{O(P)}$ does not meet the center of the 
  Wronski projection $\pi_\alpha$.
\end{lemma}

\noindent{\it Proof.}
  We show that on $Y_{O(P)}$ the equations defining the center of the projection
  $\pi_\alpha$
\[
   \sum_{|J|=k}\alpha_{|J|}x_{J}\ =\ 0\,,\qquad k=0,\dotsc,n\,,
\]
 generate the irrelevant ideal $\langle x_J\mid J\in\calJ(P)\rangle$.

 Let $I$ be the ideal of the equations for the center of projection and the
 equations defining $X_\Delta$.
 If there is only one order ideal $J$ with $|J|=k$, then $x_J\in I$,
 in particular, $x_J\in I$ when $|J|=0$.
 Suppose that we have $x_{J}\in I$ for all $J$ with $|J|<k$.
 Given two order ideals $J$ and $K$ of size $k$,  we have
 $x_J x_K=x_{J\wedge K}x_{J\vee K}\in I$ as $|J\vee K|<k$.
 Together with the equation defining the center of projection,
 this implies that if $|J|=k$, then $x_J\in I$.
 By induction on $k$, $I=\langle x_J\mid J\in\calJ(P)\rangle$.

 This argument also shows that  $s.Y_{O(P)}$ does not meet the center of
 projection.
\hfill \raisebox{-3pt}{\epsfxsize=12pt\epsffile{figures/QED.eps}}\medskip

\begin{lemma}
 If a finite poset $P$ is ranked mod $2$, then $Y_{O(P)}^{+}$ is Cox-oriented.
\end{lemma}

\noindent{\it Proof.}
Let the order polytope be defined by facet inequalities
\[
   O(P)\ =\ \big\{y\in\R^n\mid\mathcal{A} y\geq-b \big\}\,,
\]
where $\mathcal A$ is an integer $r\times n$ matrix.
By Theorem~\ref{T:orientability} it is enough to check that
the vector consisting of all ones is in the mod 2 integer column span of the
matrix $[\mathcal{A}:b]$ and the lattice $\Lambda_{\mathcal{A}}$ spanned by the
columns of $\mathcal A$ is saturated.
The integral points of $O(P)$ affinely span $\Z^P$ as $O(P)$ has a unimodular
triangulation.

Each facet of the order polytope $O(P)$ is defined by one of the following conditions:
$$\begin{array}{cll}
&y_a=0\quad &\textrm{for a minimal}\ a\in P,\\
&y_b=1\quad &\textrm{for a maximal}\ b\in P, \\
&y_a=y_b\quad &\textrm{for}\ a\  \textrm{covering}\ b\ \textrm{in}\ P.
\end{array}
$$
Fix a maximal chain $a_1<\dots<a_k$ in $P$.
The corresponding facets of $O(P)$
are
\[
 y_{a_1}\ =\ 0\,,\quad y_{a_2}-y_{a_1}\ =\ 0\,,\quad\dots,\quad 
 y_{a_k}-y_{a_{k-1}}\ =\ 0\,,
\quad y_{a_k}\ =\ 1\,,
\]
and the corresponding rows of the matrix $[\mathcal{A}:b]$ are
\[
 \left[
 \begin{array}{rrrcrrrcc}
 1&0&0&\dots&0&0&0&\dots&0\,\\
 -1&1&0&\dots&0&0&0&\dots&0\,\\
 0&-1&1&\dots&0&0&0&\dots&0\,\\
 \vdots &\vdots &\vdots &\ddots&\vdots &\vdots&\vdots&\dots&\vdots\,\\
 0&0&0&\dots&1&0&0&\dots&0\,\\
 0&0&0&\dots&-1&1&0&\dots&0\,\\
 0&0&0&\dots&0&1&0&\dots&0\,
 \end{array}\,\right|\,
 \left.
 \begin{matrix}
 0\\
 0\\
 0\\
 \vdots \\
 0\\
 0\\
 1
 \end{matrix}\right]\ .
\]
The columns of $\mathcal A$ are indexed by the elements of $P$.
Consider the linear combination of the columns of $\mathcal A$ where the
coefficient of a column corresponding to an element $a$ of $P$ is
$1-\mbox{rk}(a)$, where $\mbox{rk}(a)$ is its mod 2 rank.
This will be a vector with all components odd if $P$ has mod 2 rank 0.
If $P$ has mod 2 rank 1, then adding the vector $b$ to this combination gives
a vector with all components odd.
(Here, a minimal element has rank 0.)

Consider the submatrix of $\calA$ consisting of its rows corresponding to
minimal elements of $P$, together with one row for each non-minimal element $a$
of $P$ corresponding to some cover $a'\lessdot a$.
This submatrix has determinant $\pm 1$, which implies that
column space of the matrix $\mathcal A$ is saturated.
\hfill \raisebox{-3pt}{\epsfxsize=12pt\epsffile{figures/QED.eps}}\medskip

We noted earlier that $\ch(f)=0$ whenever $Y_{\Delta}$ is orientable
but $\RP^n$ is not. By Remark~\ref{R:orient_Y} $Y_{\Delta}$ is orientable
if there exists a vector all of whose components are odd in the integer
column span of $\calA$. For posets, this translates to: $Y_{O(P)}$ is orientable
if all the maximal chains of $P$ are odd. We obtain:

\begin{cor}
If all the maximal chains of a finite poset $P$ are odd but the number of
elements in $P$ is even then the poset $P$ is sign-balanced.
\end{cor}

This is Corollary 2.2 of \cite{St02}, where it is given a a purely
combinatorial proof.

\section{Further examples}

This theory applies to other toric varieties besides those associated to the
order polytopes of finite posets.
The hexagon of Example~\ref{Ex:Hexagon} is one instance.
We present three additional instances based on particular triangulations of
polytopes, and one infinite family that is based on the chain polytopes
of~\cite{St86b}.

\subsection{Three examples of polytopes}

\begin{ex}\label{Ex:Triangle}
  Let $\Delta$ be the convex hull of the points $(0,0)$, $(0,3)$, and 
  $(3,0)$, a triangle.  
  Then $T_\Delta$ is a Veronese embedding of $\R\P^2$ and $Y^+_\Delta$ is the 
  2-sphere, and so it is orientable.
  The triangle has a regular unimodular balanced triangulation with
  signature 3 illustrated in Figure~\ref{F:triangle} below. 
  This is induced by a weight function whose values are 0 at the center, 3 at
  each of the three vertices, and 1 at the remaining six points.
  A Wronski polynomial with constant weight 1 on members of the family
  $s.Y_\Delta$ has the form 
 \begin{equation}\label{E:Triang}
   a(s^3+xy+s^3x^3+s^3y^3) + bs(x+x^2y+y^2) + cs(y+xy^2+x^2)\ =\ 0\,.
 \end{equation}
 The center of projection does not meet $s.Y_\Delta$, for any $s\neq 0$.
%
%
 Thus any two polynomial equations of the form~\eqref{E:Triang} will
 have at least 3 common real solutions.
 Figure~\ref{F:triangle} also shows two curves given by equations of the
 form~\eqref{E:Triang} with $s=1$ and coefficients $\Blue{(4,-11,4)}$ and 
 $\Red{(-13, -1, 24)}$.  These meet in 9 points, giving 9 solutions to the
 system. 
 \begin{figure}[htb]
 \[
   \begin{picture}(120,108)
     \put(12,10){\epsfxsize=1.3in\epsffile{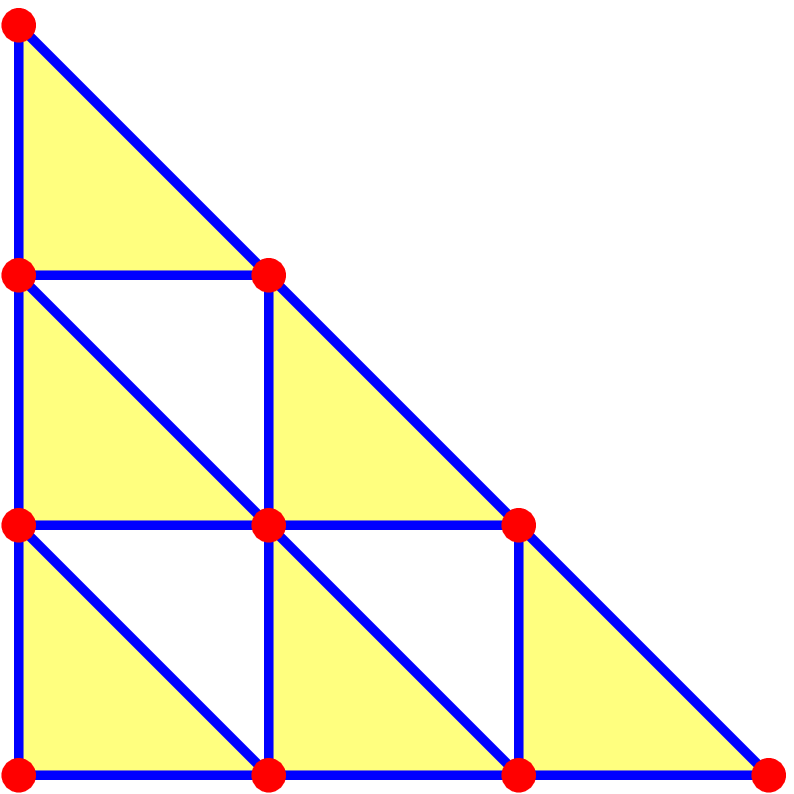}}
     \put(0,98){$y^3$}
     \put(0,70){$y^2$} \put(47,76){$xy^2$}
     \put(2,40){$y$}\put(29, 33){$xy$}   \put(77,46){$x^2y$}
     \put(5, 0){$1$}\put(40,0){$x$}\put(70,0){$x^2$}\put(100,0){$x^3$}
     \end{picture}
    \qquad \qquad 
   \epsfxsize=1.5in\epsffile{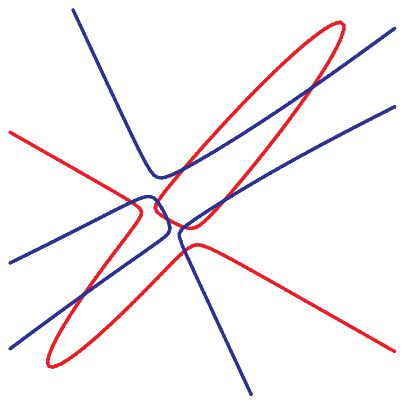} 
 \]
 \caption{Cubic system}
 \label{F:triangle}
 \end{figure}
 We computed ten million instances of the Wronski system on $Y_\Delta$, where
 the coefficients $a,b,c$ were integers chosen uniformly from the interval
 $[-1000,1000]$. 
 Of these, 9,\,976,\,701 ($99.8\%$)  had 3 real solutions and 23,\,299 had 9
 real solutions. 
 Computing 500,\,000 instances of the Wronski system~\eqref{E:Triang} with
 $s\in(0,1)$, we found 414,\,592 with 3 real solutions, 85,\,408 with 9 real
 solutions, and did not find any with either 5 or 7 real solutions. 

\end{ex}

\begin{ex}\label{Ex:Cube_Unimodular}
 Let $\Delta$ be the unit cube.
 Then $Y_\Delta=(S^1)^3$ and so it is oriented.
 Consider the regular unimodular balanced triangulation of the unit
 cube $[0,1]^3$ illustrated on the left in Figure~\ref{F:cube}.
 It has signature $4-2=2$ and is 
 given by a weight function taking values 3 at $(0,0,0)$ and $(1,1,1)$, 
 0 at $(1,0,0)$ and $(0,1,1)$, and 1 at the remaining vertices.
 The corresponding Wronski polynomials on $s.Y_\Delta$ 
 have the form
\[
  a(s^3+yz) + b(x+s^3xyz) + cs(y+xz) + ds(z+xy)\,.
\]
 The family $s.Y_\Delta$ meets the centre of projection only when $s^3=\pm1$.
 These points for $s=1$ are  
\[
   (x,y,z)\ \in\ \{(1,1,-1), (1,-1,1), (1,i,i), (1,-i,-i)\}\,.
\]
 Thus $Y_\Delta$ meets the center of projection in 2 real and 2 complex points.
 Theorem~\ref{T:computedegree} implies that for $s\in(0,1)$, there will be
 at least 2 real solutions, and we have computed such systems with 2, 4, and 6
 real solutions. 
\begin{figure}[htb]
 \[
   \begin{picture}(150,145)(0,-20)
    \put(25,10){\epsfxsize=1.5in\epsffile{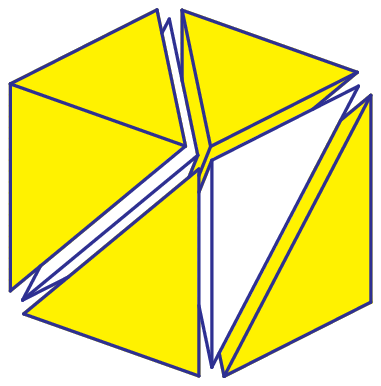}}
                  \put(65,125){$xz$}
    \put(0,95){$xyz$}          \put(135,100){$z$}
                  \put(90,67){$yz$}
    \put(10,22){$xy$}          \put(138,25){$1$}
                  \put(80,0){$y$}
    \put(50,-20){Unimodular}
   \end{picture}
    \qquad
   \begin{picture}(147,145)(5,-18)
    \put(25,10){\epsfxsize=1.4in\epsffile{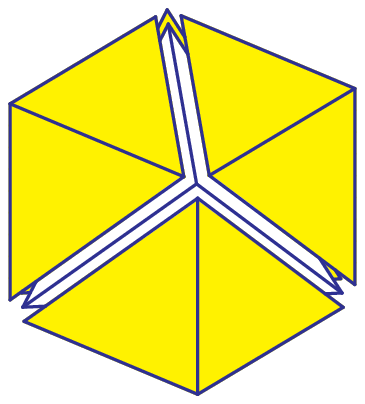}}
                  \put(65,129){$xz$}
    \put(3,95){$xyz$}          \put(130,100){$z$}
                  \put(83,48){$yz$}
    \put(14,27){$xy$}          \put(129,33){$1$}
                  \put(78,2){$y$}
    \put(43,-20){Non-unimodular}
   \end{picture}
\]
 \caption{Triangulations of Cubes}
 \label{F:cube}
 \end{figure}
\end{ex}

\begin{ex}

 On the right of Figure~\ref{F:cube} is the regular triangulation of the cube 
 given by the lifting function that takes values 0 at $(0,0,0)$, $(1,1,0)$,
 $(1,0,1)$, and $(0,1,1)$, and 1 at the remaining four vertices.
 It is balanced with 4 unimodular simplices of the same color and one with
 normalized volume 2 of the opposite color, and thus has signature 4.
 The corresponding Wronski polynomials on $s.Y_\Delta$ have the form
\[
  a(1+sxyz) + b(sx+yz) + c(sy+xz) + d(sz+xy)\,.
\]
 The variety $s.Y_\Delta$ meets the center of projection only when $s^4=1$.
 When $s=1$, there are four points of intersection
\[
   \{ (x,y,z)\mid x,y,z\in\{\pm1\}, xyz=-1\}\,.
\]

 Since the sign imbalance is 4, Theorem~\ref{T:computedegree} implies that 
 for $s\in(0,1)$, three polynomials of this form
 will have at least 4 real solutions.
 Computing 500,000 instances, we found 453,811 with 4 solutions and 46,189
 with 6 solutions.  
%
\end{ex}
%
%
\subsection{Systems from chain polytopes}\label{S:chainPolytope}

 Let $P$ be a poset with $n$ elements.
 Stanley~\cite{St86b} defined the {\it chain polytope} $C(P)$ to be the set of
 points $y$ in the unit cube $[0,1]^P$ such that 
\[
  y_a+y_b+\dotsb+y_c\leq 1\qquad\textrm{whenever}\qquad 
   a<b<\dotsb<c\quad \textrm{is a chain in }P\,.
\]

 This polytope is intimately related to the order polytope $O(P)$ of
 Section~\ref{S:OrderPolytope}. 
 It has no interior lattice points but its vertices are the
 characteristic functions of the antichains of $P$, 
 and the bijection between (upper) order ideals $J$ and antichains $A$ given by
 \begin{eqnarray*}
  J&\longmapsto& \textrm{minimal elements in }J\\
    A&\longmapsto&\langle A\rangle\ :=\ 
   \{b\in P\mid a\leq b, \textrm{ for some }a\in A\}
 \end{eqnarray*}
 extends to a bijection $\varphi$ between the polytopes.
 Let $y\in O(P)$ be a point in $[0,1]^P$ with $y_a\leq y_b$ whenever  
 $a\leq b$ in $P$. 
 For $a\in P$, define
 \begin{equation}\label{E:varphi}
   \varphi(y)_a\ =\ \min\{y_a-y_b\mid a \textrm{ covers $b$ in }P\}\,.
 \end{equation}
This is piecewise linear on the simplices of the canonical 
triangulation of $O(P)$ and it extends the bijection given above. 
This induces the {\it canonical triangulation} of the chain polytope, which is
unimodular, balanced, and has the same 
signature as the canonical triangulation of the order polytope.
It is regular, by Lemma~\ref{L:ChainRegular}.

 Let $\calA(P)$ denote the set of antichains of $P$.
 Let ${\rm rk}(A)$ be the number of elements in the (upper)
 order ideal generated by the antichain $A$.
 For an antichain $A$, let $t^A:=\prod_{a\in A} t_a$ be the
 monomial in $\R[t_p\mid p\in P]$ whose exponent vector is the vertex of $C(P)$
 corresponding to the antichain $A$.
 Fix a system of weights $\{\alpha_A\in\R^\times\mid A\in\calA(P)\}$.
 Given a coefficent vector  $c=(c_0,c_1,\dots,c_n)\in(\R^\times)^{n+1}$,
 set
  \begin{equation}\label{E:chain_poly}
   F_c(t)\ :=\ \sum_{A\in\calA(P)} c_{{\rm rk}(A)} \alpha_A t^A\,.
  \end{equation}
A system of $n$ such polynomials for a fixed choice $\alpha$ of weights is a
Wronski polynomial system for the canonical triangulation of the chain polytope
of $P$ with weight $\alpha$.

\begin{thm}\label{T:Chain_Poly}
 Suppose that a finite poset $P$ is ranked mod $2$.
 For any choice of weights, a Wronski polynomial system for the canonical
 triangulation of the chain polytope of $P$ with weight $\alpha$ will have
 at least $\sigma(P)$ real solutions.
\end{thm}

 In Section~\ref{S:gaps}, we consider such systems when $P$ is a
 incomparable union of chains.

\begin{ex}
 Despite the similarities between our results for the chain and order polytopes,
 the polytopes $O(P)$ and $C(P)$ are not isomorphic, in general.
 For example, if $B_n$ is the boolean poset $\{0,1\}^n$, then 
 $B_n$ has $2^n$ elements and $n!$ maximal chains.  It also has 
 unique maximal and minimal elements and exactly $4\cdot 3^{n-2}$ covers.
 Thus, for $n\geq 2$, 
  \[
     \textrm{$O(B_n)$ has $2{+}4{\cdot}3^{n-2}$ facets, \ while 
     $C(B_n)$ has $2^n{+}n!$ facets.}
  \]
%
 In particular, $O(B_4)$ has 38 facets while $C(B_4)$ has 40.
 When $n=10$, these numbers of facets are $26,246$ and $3,629,824$,
 respectively. 
\end{ex}

The following four lemmas reduce Theorem~\ref{T:Chain_Poly} to
Theorem~\ref{T:computedegree}.   

\begin{lemma}\label{L:ChainRegular}
  The canonical triangulation of the chain polytope is regular.
\end{lemma}

 We defer the proof until the end of this section.
 There, we will show that the lifting function $\omega(A):=3^{{\rm rk}(A)}$
 induces the canonical triangulation of the chain polytope.

\begin{lemma}
  The canonical triangulation of the chain polytope $C(P)$ is balanced. 
  Its signature is $\sigma(P)$, the sign imbalance of $P$.
\end{lemma}

\noindent{\it Proof. }
 The map $\varphi$ between simplices in the canonical
 triangulations of the chain and order polytopes shows that the two
 triangulations are combinatorially equivalent.
 The statement then follows from Lemma~\ref{L:Order_balanced}.
\hfill \raisebox{-3pt}{\epsfxsize=12pt\epsffile{figures/QED.eps}}\medskip

 The lifting function $\omega$ which induces the canonical triangulation of
 $C(P)$ give a $\R^\times$-action on the projective space
 $\R\P^{C(P)}$ in which $Y_{C(P)}$ lives.
 Let $s.Y_{C(P)}$ be the associated toric deformation of $Y_{C(P)}$.

\begin{lemma}
 For any choice of weights $\{\alpha_A\mid A\in\calA(P)\}$, the toric
 degeneration of\/$Y_{C(P)}$ induced by the lifting function 
 $\omega(A):=3^{{\rm rk}(A)}$  does not meet the center of the projection
 $\pi_\alpha$ defining the Wronski polynomial system.
\end{lemma}

\noindent{\it Proof.}
 The center of the projection $\pi_\alpha$ is defined by the equations
\[ 
   \sum_{{\rm rk}(A)=k} \alpha_A t^A\ =\ 0\qquad
    \mbox{for } k=0,1,\dotsc,n=|P|\,.
\]
 As in the proof of Lemma~\ref{L:center_proj_orderpoly}, it suffices to show
 that if $A$ and $B$ are two antichains with the same rank $k>0$, then there is
 a relation in the ideal of $Y_{C(P)}$ of the form
 \begin{equation}\label{E:relation_chain}
   x_A x_B\ -\ x_C x_D\,,
 \end{equation}
 where $C,D$ are antichains with ${\mbox rk}(C)<k$.

 Let $A$ and $B$ be two antichains, each of rank $k$.
 Let $D$ be the antichain of minimal elements in 
 $\langle A\rangle\cup\langle B\rangle$.
 Then $D\subset A\cup B$.
 Set $C := [(A\cup B) - D] \cup (A\cap B)$.
 This is a subset (possibly proper) of the minimal elements of
 $\langle A\rangle\cap\langle B\rangle$.
 Since $A\neq B$, we have $|\langle A\rangle\cap\langle B\rangle|<k$, and so 
 $C$ is an antichain with rank less than $k$.
 Finally, the multiset inequality $A\cup B=C\cup D$ implies the
 relation~\eqref{E:relation_chain}, which completes the proof.
\hfill \raisebox{-3pt}{\epsfxsize=12pt\epsffile{figures/QED.eps}}\medskip

\begin{rem}
  The sets $C$ and $D$ may be constructed from any two incomparable antichains
  $A$ and $B$ of $P$.
  From the construction, we have
  ${\rm rk}(A)+{\rm rk}(B)\geq{\rm rk}(C)+{\rm rk}(D)$,
  with equality only when $C$ is the set of minimal elements of
  $\langle A\rangle\cap\langle B\rangle$.

 Despite the similarity of the equations for the order
 polytope~\eqref{E:toric_from_lattice} to those in the proof
 above~\eqref{E:relation_chain} for the chain polytope,
 the equations for the chain polytope do not come from a distributive lattice. 
\end{rem}

 \begin{lemma}
 If a finite poset $P$ is ranked mod $2$, then $Y_{C(P)}^{+}$ is Cox-oriented.
 \end{lemma}

 \noindent{\it Proof.}
  The chain polytope contains the standard basis of $\R^P$ as some of its
  vertices and it contains the origin.
  Thus its integral points affinely span $\Z^P$.
  The facet inequalities $\calA x\geq b$ for $C(P)$
  come in two forms
  \begin{eqnarray*}
    f(a)\geq 0&\textrm{for}& a\in P, \ \textrm{and}\\
    f(a)+f(b)+\dotsb+f(c)\leq 1&\textrm{for}&
     a<b<\dotsb<c\textrm{ a maximal chain in }P
  \end{eqnarray*}
  The first collection of inequalities ensures that the columns of the matrix
  $\calA$ span an $n$-dimensional saturated sublattice of $\Z^{n+c}$, where $c$
  is the number of maximal chains.
  If $P$ has mod 2 rank 0, then the sum of the columns of $\calA$ is a vector in
  $Z^{n+c}$ with every component odd. 
  If $P$ has mod 2 rank 1, then we add $b$ to the sum of columns of $A$ 
  gives a vector in $Z^{n+c}$ with every component odd.

  The lemma follows by Theorem~\ref{T:orientability}.
  \hfill \raisebox{-3pt}{\epsfxsize=12pt\epsffile{figures/QED.eps}}\medskip

\noindent{\it Proof of Lemma~$\ref{L:ChainRegular}$.}
 For an antichain $A$ of $P$ generating an order ideal with $k$ elements,
 set $\omega(A):=3^k-1$.
 We show that $\omega$ induces the canonical triangulation.

 Suppose that $P$ has $n$ elements and fix a linear extension 
 $\pi\colon P\to[n]$. 
 Let $\Delta_\pi$ be the simplex in the canonical triangulation corresponding to
 this linear extension.
 For each $k=0,\dotsc,n$, let $I_k$ be the order ideal
 $\pi^{-1}\{n{+}1{-}k,\dotsc,n{-}1,n\}$ and let $A_k$ be the antichain of
 minimal elements of $I_k$. 
 Then the vertices of $\Delta_\pi$ are $m_k=\sum_{x\in A_k}e_x$, for
 $k=0,\dotsc,n$, where
 $\{e_x\mid x\in P\}$ is the standard basis for $\R^P$. 

 Let $\Lambda$ be the unique linear function satisfying 
 $\Lambda(m_k)+\omega(m_k)=0$. 
 We will show that if $m$ is a vertex of the chain polytope but not of
 $\Delta_\pi$, then $\Lambda(m)+\omega(m)<0$, which will prove the proposition.
 This requires a more precise description of $\Lambda$.
 Write $p\lessdot q$ if $q$ covers $p$ in $P$.
 For $x,y\in P$, define the function $\beta_{x,y}\in\{0,1\}$ recursively as
 follows. 
\begin{enumerate}
 \item $\beta_{x,y}=0$ if $x\not\leq y$, 
 \item $\beta_{y,y}=1$, and 
 \item $\beta_{x,y}=\sum\{\beta_{z,y}\mid x\lessdot z,\ \textrm{and}\ 
                           z\in A_{\pi(x)+1}\}$.
\end{enumerate}

 These functions $\beta_{x,y}$ have the following elementary and obvious
 properties. 

\begin{lemma}\label{L:alpha}
 Let $\pi\colon P\to [n]$ be a linear extension and define $\alpha$ as above.
 Then

 \begin{enumerate}
  \item For any $y\in P$ and antichain $A$, 
        ${\displaystyle \sum_{x\in A} \beta_{x,y}\leq 1}$.
  \item If $\pi(y)\geq j$, then
        ${\displaystyle 1=\sum_{x\in A_j} \beta_{x,y}}$.
 \end{enumerate}
\end{lemma}

 Set $f(k):=3^{k-1}-3^k$ and note that if
\[ 
   \Lambda(e_x)\ :=\ \sum_y \beta_{x,y}f(n-\pi(y)+1)\,,
\]
 then $\Lambda(m_k)=1-3^k$.
 Indeed, 
\begin{eqnarray*}
  \Lambda(m_k)&=&\sum_{x\in A_k}\Lambda(e_x)\ =\ 
                 \sum_{x\in A_k} \sum_y\beta_{x,y}f(n-\pi(y)+1)\\
              &=&\sum_{y\in I_k} f(n-\pi(y)+1)\sum_{x\in A_k} \beta_{x,y}\\
              &=&\sum_{y\in I_k} f(n-\pi(y)+1)\ =\ \sum_{i=0}^k f(k)\ =\ 1-3^k\,.
\end{eqnarray*}

 Suppose now that $m$ is a vertex of $C(P)$, but $m\not\in\Delta_\pi$.
 Let $A$ be the antichain corresponding to $m$ and let $z$ be the least
 element of $A$ under the linear extension $\pi$.
 Set $k:=n+1-\pi(z)$.
 Then $A\subset I_k$ but $A$ does not generate $I_k$ (for otherwise
 $m=m_k\in\Delta_\pi$), and so $\omega(m)\leq 3^{k-1}-1$, as the order ideal
 generated by $A$ has at most $k-1$ elements.
 However, $e_z$ occurs in $m$, so 
 $\Lambda(m)\leq f(n+1-\pi(z))=f(k)=3^{k-1}-3^k$.
 But then $\Lambda(m)+\omega(m)<0$, as claimed.
\hfill \raisebox{-3pt}{\epsfxsize=12pt\epsffile{figures/QED.eps}}\medskip

\section{Lower bounds from sagbi degeneration}

The Grassmannian has a flat deformation to the toric variety of an order
polytope induced by a canonical subalgebra or sagbi basis of its
coordinate ring.
We use this sagbi deformation to compute the characteristic of the
real Wronski map and recover the result of Eremenko and
Gabrielov~\cite{EG01b} which motivated our work. 
More generally, we compute the characteristic of Wronski
projections for many projective varieties whose coordinate rings are
algebras with a straightening law on a distributive lattice.

We review some definitions from \cite{EG01b}. 
Let $f_1(z),\dots,f_p(z)$ be real polynomials in one variable $z$, each of
degree at most $m{+}p{-}1$.
Their Wronski determinant is 
\[ 
 \begin{matrix}
     W(f_1(z),\dots,f_p(z))\ :=\ 
  \end{matrix}
   \left|
  \begin{array}{ccc}
   f_1(z) &\dotsb  &f_p(z)  \\
   f_1'(z)&\dotsb  &f'_p(z) \\
   \vdots &       &\vdots \\ 
   f_1^{(p-1)}(z)&\dotsb&f^{(p-1)}_p(z)
  \end{array}\right|\ .
\]
This Wronskian has degree at most $mp$, and,
up to a scalar factor, it depends only upon the linear span of the
polynomials $f_1(z), f_2(z),\dotsc,f_p(z)$.
If we identify a polynomial of degree $m{+}p{-}1$ as a linear form on $\R^{m+p}$,
then $p$ linearly independent polynomials cut out a $m$-plane.
Thus the Wronski determinant induces the {\it Wronski map}
\[
   W\ \colon\ G(m,p)\ \longrightarrow\ \P^{mp}\ ,
\]
where $G(m,p)$ is the Grassmannian of $m$-planes in $\R^{m+p}$, and $\P^{mp}$ is
the space of polynomials of degree $mp$, modulo scalars.

Consider this in more detail.
Represent a polynomial $f(z)$ by the column vector $f$ of its coefficients.
Set
\[
   \gamma\ :=\ [1,z,z^2,\dotsc,z^{m+p-1}]\,,
\]
and define $K(z)$ to be the matrix with rows 
$\gamma(z), \gamma'(z),\dotsc, \gamma^{(m-1)}(z)$.
Then
\[
   \left[\begin{matrix} 
      f_1(z)& \dotsb & f_p(z)\\
      \vdots&\ddots & \vdots \\
   f^{(m-1)}_1(z)&\dotsb&f^{(m-1)}_p(z)
    \end{matrix}\right]
    \ =\ 
   K(z)\cdot [f_1,\dotsc,f_p]\,.
\]
Expanding the determinant using the Cauchy-Binet formula gives
\[
    W(f_1(z),\dots,f_p(z))\ =\
    \sum_J k_J(z) x_J(f_1,\dotsc,f_p)\,,
\]
where the summation is over all sequences 
$J: 1\leq j_1<\dotsb<j_p\leq m{+}p$,
$x_J(f_1,\dotsc,f_p)$ is the determinant of the $p\times p$ submatrix of 
$[f_1,\dotsc,f_p]$ formed by the rows in $J$, and $(-1)^{mp-|J|}k_J(z)$ is 
the determinant of the complementary rows of $K(z)$.
These functions $x_J(f_1,\dotsc,f_p)$ are the {\it Pl\"ucker coordinates} of the
$m$-plane cut out by $f_1,\dotsc,f_p$.
They define a projective embedding of $G(m,p)$ into Pl\"ucker space,
$\P^N$, where $N=\binom{m{+}p}{m}-1$.

If $|J|:=\sum j_i-i$, then $k_J(z)=z^{mp-|J|}k_J(1)$.
Moreover, $(-1)^{mp-|J|}k_J(1)>0$ for all $J$ (Equation (5.5)
of~\cite{So03}\footnote{There is a misprint in the cited paper at this point,
       $F_j(0)$ should be $F_i(1)$.}).
If we write $\alpha_J:=k_J(1)$, then, in the Pl\"ucker coordinates for $G(m,p)$
and the basis of coefficients for polynomials in $\P^{mp}$, the Wronski map is
\[
   W(x_J\mid J\in \calC_{m,p})\ =\ 
   \sum_{j=0}^{mp} z^{mp-j}\ 
   \sum_{|J|=j} \alpha_J x_J\,,
\]
where $\calC_{m,p}$ is the set of indices of Pl\"ucker coordinates.
We recognize this as the restriction of a linear projection $\pi_\alpha$ 
on $\P^{mp}$ to the Grassmannian $G(m,p)$.

\begin{figure}[htb]
  \centerline{
  \scalebox{1.0}{
  \input{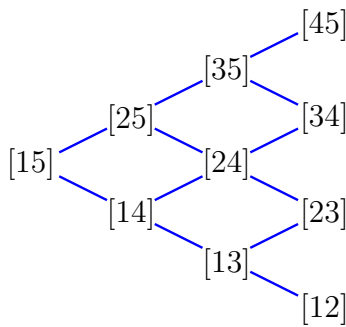}
             }
           }
  \caption{The distributive lattice $\calC_{2,3}$.}\label{F:C23}
 \end{figure}

This is a Wronski projection, in the sense of~\ref{S:char}.
Indeed, indices $\calC_{m,p}$ of Pl\"ucker coordinates are partially ordered 
by componentwise comparison.
\[
  J\ =\ [j_1,\dots,j_p]\ \leq\ K\ =[k_1,\dots,k_p]\qquad
   \iff\quad j_i\leq k_i \ \ \text{for}\ \ i=1,\ \dots,p\,.
\]
This poset $\calC_{m,p}$ is the lattice of order ideals of the poset 
$[m]\times[p]$ of two chains of lengths $m$ and $p$ and the rank of $J$ is $|J|$.
Figure~\ref{F:C23} shows $\calC_{3,2}$.

Following~\cite{EG01b} we define the characteristic of the Wronski map $W$.
This map sends the subset $X$ of $G(m,p)$ where 
$x_{1,2,\dotsc,p}\neq 0$ to the subset $Y$ of $\P^{mp}$ of monic polynomials of
degree $mp$, and the complement of $X$ to the complement of $Y$.
Both $X$ and $Y$ are orientable as they are identified with $\R^{mp}$. 
Let $\ch(W)$ be the absolute value of the topological
degree of $W|_X$.

An equivalent definition of $\ch(W)$ similar to Definition~\ref{D:char} also
appears in~\cite{EG01b}. 
Lift the Grassmannian and the projection to the double covers $S^N$ and $S^{mp}$
of $\RP^N$ and $\RP^{mp}$. 
The pullback of $G(m,m)$ is the upper Grassmannian $G^{+}(m,m)$ of all
oriented $m$-planes in $\R^{m+p}$. 
Let $W^{+}\colon G^{+}(m,m)\to S^{mn}$ be the pullback of $W$. 
Since $G^{+}(m,m)$ is orientable~\cite{EG01b}, $\ch(W)$ is
well-defined as in Definition~\ref{D:char}.

\begin{thm}\cite[Theorem 2]{EG01b},
 The characteristic of the real Wronski map $W$
 is equal to the sign-imbalance of $\calC_{m,p}$.
\end{thm}

\begin{rem}
 White~\cite{Wh01} computed this sign-imbalance, showing that
 $\sigma(\calC_{m,p})=0$ unless $m+p$ is odd, and then it equals
\[
   \frac{1!2!\dotsb(p{-}1)!(m{-}1)!(m{-}2)!\dotsb(m{-}p{+}1)!(\frac{mp}{2})!}
   {(m{-}p{+}2)!(m{-}p{+}4)!\dotsb(m{+}p{-}2)!\left(\frac{m-p+1}{2}\right)!%
     \left(\frac{m-p+3}{2}\right)!\dotsb\left(\frac{m+p-1}{2}\right)!}\ .
\]
\end{rem}

\noindent{\it Proof.}
The Pl\"ucker ideal $I$ of the Grassmannian $G(m,p)$ has a quadratic Gr\"obner
basis whose elements are indexed by incomparable pairs $J,K$ in $\calC_{m,p}$
and have the form
 \begin{equation}\label{Eq:PlueckerEqs}
   x_J x_K\ -\ x_{J\wedge K}x_{J\vee K}\ +\ \mbox{\rm other terms}\,,
 \end{equation}
where the other terms have the form $ax_Lx_M$ with 
$L\leq J,K \leq M$~\cite[Chapter 3]{Sturmfels_Alg}.
The term order here is degree reverse lexicographic on 
$\C[x_J\mid J\in\calC_{m,p}]$ where the variables are first linearly ordered by
the ordinary lexicographic order on their indices.

Any lifting function $\omega\colon\calC(m,p)\to \Z$ defines a $\R^\times$-action
on Pl\"ucker space by $s.x_J=s^{-\omega(J)}x_J$. 
Restricting this action to the Grassmannian gives a family $s.G(m,p)$ whose
scheme theoretic limit as $s\to 0$ is cut out by the initial ideal 
$\ini_\omega I$.
There is a lifting function $\omega$ so that $\ini_\omega I$ is the toric ideal
\[
   x_J x_K\ -\ x_{J\wedge K}x_{J\vee K}\,,\qquad\mbox{for $J,K$ incomparable}
\]
of the distributive lattice $\calC_{m,p}$~\cite[Theorem 11.4]{Sturmfels_GBCP}.
We call the corresponding family $s.G(m,p)$ the {\it sagbi} deformation of the
Grassmannian, which deforms it into the toric variety $Y(m,p)$ of the distributive
lattice $\calC_{m,p}$.
As in the proof of Lemma~\ref{L:center_proj_orderpoly}, the form of the
equations~\eqref{Eq:PlueckerEqs} for the Grassmannian imply that 
the sagbi degeneration does not meet the center of the projection $\pi_\alpha$.

For $s\in(0,1]$, let $W_s$ be the restriction of the Wronski $\pi_\alpha$ to
$s.G(m,p)$. 
Then the characteristic $\ch(W)$ coincides with $\ch(W_s)$ for $s\in(0,1]$. 
By Lemma~\ref{L:pos_comp}, there exists a regular value $z$ of $W$ all of whose
preimages in $Y(m,p)$ have all components positive. 
By the implicit function theorem, every preimage of $z$ in
$s.G(m,p)$ for $s$ suficiently small have the same property.
Thus $\ch(W)$ equals to the
degree of the restriction of $W$ to $s.G_>(m,p)$, the intersection of the
Grassmannian with the positive orthant. 

Let $Y_>(m,p)$ be the positive part of the toric variety $Y(m,p)$.
Then $Y_>(m,p)$ has coordinates and orientation defined by the projection to the
coordinate plane corresponding to some (fixed) maximal chain in $\calC_{m,m+p}$. 
By the implicit function theorem, the same is true for $s.G_>(m,p)$ when $s$ is
sufficently small. 
For each preimage of $z$ in $Y_>(m,p)$ there is a nearby preimage of $z$ in 
$s.G_>(m,p)$. 
Hence the projections of these preimages to the coordinate plane
are nearby, and the signs of $\det W_s$ and $\det W_0$ coincide. 
This proves that $\ch(W)=\ch\bigl(W_0|_{Y_>(m,p)}\bigr)$. 
Finally, by Corollary~\ref{C:Y_positive_char}, 
$\ch (W_0|_{Y_>(m,p)})$ is equal to the sign-imbalance
of $\calC_{m,p}$, which completes the proof.
\hfill \raisebox{-3pt}{\epsfxsize=12pt\epsffile{figures/QED.eps}}\medskip

We only used that the characteristic of the Wronski map was defined and that 
$G(m,p)$ has equations of the form~\eqref{Eq:PlueckerEqs}, for some distributive
lattice $D$. 
The projective coordinate ring of a variety $Y$ with such equations is an
algebra with straightening law on the distributive lattice $D$~\cite{CEP82}, and
this has the geometric consequence that $Y$ admits a flat degeneration to the
toric variety $Y_D$ of the distributive lattice.
There are many examples of such varieties, besides the Grassmannian.
These include Schubert varieties of the Grassmannian, the classical flag
variety, and the Drinfel'd compactification of the space of curves on the
Grassmannian~\cite{SS01}, as well as products of such spaces.

Such a variety $Y$ has projective coordinates $\{x_J\mid J\in D\}$.
Given a set $\alpha$ of weights, a Wronski map for the lattice $D$ is a linear
projection $\pi_\alpha$ of the form 
\[
   \pi_\alpha\ \colon\ (x_J\mid J\in D)\ \longmapsto\ 
    \biggl(\sum_{|J|=j} \alpha_J x_J\ \mid\, j=0,\dotsc,\mbox{rank}(D)\biggr)
\]
We say that $\pi_\alpha$ has constant sign if the sign of $\alpha_J$ depends
only upon $|J|$.
Let $\hat{0}$ be the unique minimal element in $D$.
For $\calC_{m,p}$, this is $(1,2,\dotsc,p)$.

\begin{thm}\label{T:sagbi}
 Let $Y$ be a projective variety whose coordinate ring is an algebra
 with straightening law on a distributive lattice $D$
 and let $\pi_\alpha$ be a Wronski projection for this lattice with
 constant sign.
 If either $Y^+$ or $Y\cap\{x\mid x_{\hat{0}}\neq 0\}$ are oriented, then 
 the characteristic of the Wronski projection on $Y$ is equal to the
 sign-imbalance of $D$.
\end{thm}

This result for Schubert varieties, the Drinfel'd compactification, and
products of such varieties was communicated to us by Eremenko and
Gabrielov, to whom it should be accredited.

%
\section{Incomparable chains and factoring polynomials}\label{S:gaps}

We give a different proof of Theorems~\ref{T:order_poly_degree}
and~\ref{T:Chain_Poly}, when the poset $P$ is a disjoint union of chains of
lengths $a_1,a_2,\dotsc,a_k$, and the weights $\alpha$ are constant. 
Our method will be to show that solutions to a general Wronski polynomial system
for the chain polytope of $P$ are certain factorizations of a
particular univariate polynomial.
This reformulation shows that there are certain numbers of real solutions to
these systems that are forbidden to occur, which is a new
phenomenon, but which we seem to have observed in Example~\ref{Ex:Triangle}.
It also proves the sharpness of the lower bound of
Theorems~\ref{T:order_poly_degree} and~\ref{T:Chain_Poly} for these posets, 
and shows that the conclusion holds even when the hypotheses of those theorems do
not, as $Y^+_{C(P)}$ is not orientable if the $a_i$ do not all have the same parity.
This analysis extends to posets which are the incomparable unions of other posets. 

Let $P$ be the incomparable union of chains of lengths $a_1,a_2,\dotsc,a_k$.
For each $i=1,\dotsc,k$, set $x_{i,0}:=1$ and let
\[
  x_{i,1}\ >\ x_{i,2}\ >\ \dotsb\ >\ x_{i,a_k}
\]
be indeterminates which we identify with the elements in the $i$th chain,
ordered as indicated.
Observe that the upper order ideal generated by $x_{i,j}$ has $j$ elements.
Antichains of $P$ correspond to monomials
\[
  x_{1,i_1} x_{2,i_2} \dotsc x_{k,a_k}\,,
\]
and the order ideal generated by this antichain has $i_1+i_2+\dotsb+i_k$
elements. 

A Wronski polynomial with constant weight 1 for the canonical
triangulation of the chain polytope $C(P)$ has the form
\begin{eqnarray}
 F&=& \sum_{i_1,\dotsc,i_k}  \label{EE:Wronski}
          c_{i_1+\dotsb+i_k}  x_{1,i_1} x_{2,i_2} \dotsc x_{k,a_k}\\
    &=& \sum_{j=0}^{a_1+\dotsb+a_k} c_t\Bigl(
    \sum_{\substack{i_1,\dotsc,i_k\\ i_1+\dotsb+i_k=j}}
        x_{1,i_1} x_{2,i_2} \dotsc x_{k,a_k}\Bigr)\ .\nonumber
\end{eqnarray}
A general system of such Wronski polynomials is equivalent to one of
the form 
 \begin{equation}\label{E:triangular}
  \sum_{\substack{i_1,\dotsc,i_k\\ i_1+\dotsb+i_k=j}}
        x_{1,i_1} x_{2,i_2} \dotsc x_{k,i_k}\ =\ b_j\qquad
   \mbox{\rm for }j=1,2,\dotsc,a_1+\dotsb+a_k\,.
 \end{equation}
Suppose that we have a solution to~\eqref{E:triangular}.
For each $i=1,\dotsc,k$, define the univariate polynomial
\[
  f_i(z)\ :=\ 1 + \sum_{j=1}^{a_i} x_{i,j} z^j\,.
\]
Then we clearly have
 \begin{equation}\label{E:factor}
  f_1(z) f_2(z)\dotsb f_k(z)\ =\ 1 + \sum_{j=1}^{a_1+\dotsb+a_k} b_j
  z^j\ =\ f(z)\,.
 \end{equation}
Similarly, any such factorization of $f(z)$ where $\deg(f_i)=a_i$
gives a solution to~\eqref{E:triangular}, and hence to our original
system. 
We have proven the following theorem.

\begin{thm}\label{T:factorization}
  The solutions to a general Wronski system with constant weights for the chain
  polytope of the incomparable union of chains of lengths $a_1,\dotsc,a_k$ 
  are the factorizations of a univariate polynomial $f$ of degree
  $a_1+\dotsb+a_k$ into polynomials $f_1,\dotsc,f_k$, where $f_i$ has degree $a_i$.
\end{thm}

\begin{rem}
 For each variable $x_{i,j}$ above, set
\[
  \varphi(x_{i,j})\ :=\ \prod_{j\leq l \leq a_i} x_{i,l}\,.
\]
 If we apply $\varphi$ to a Wronski polynomial $F$~\eqref{EE:Wronski} of
 the chain polytope of $P$, we obtain a Wronski polynomial for the canonical
 triangulation of the order polytope of $P$.
 In this way, Wronski systems for the order polytope and chain polytope of $P$
 are equivalent, and thus the results of this section also hold for the order
 polytope of $P$.
\end{rem}

We investigate the consequences of Theorem~\ref{T:factorization}.
A factorization
 \begin{equation}\label{E:New_factor}
   f_1(z) f_2(z)\dotsb f_k(z)\ =\ f(z) 
 \end{equation}
where $f_i$ is a complex polynomial of degree $a_i$ for $i=1,\dotsc,k$ and
$f(z)$ has degree $a_1+\dotsb+a_k$ and distinct roots, is a distribution
of the roots of $f$ between the polynomials $f_1,\dotsc,f_k$, with $f_i$
receiving $a_i$ roots.
Thus the number of such factorizations is the multinomial coefficient
 \begin{equation}\label{E:multinomialC}
    \binom{a_1+\dotsb+a_k}{a_1,a_2,\dotsc,a_k}\ =\ 
     \frac{(a_1+\dotsb+a_k)!}{a_1! a_2!\dotsb a_k!}\,,
 \end{equation}
which is also the number of linear extensions of $P$.
Indeed, the positions taken by the elements from each chain in a linear
extension of $P$ give a distribution of $a_1+\dotsb+a_k$ positions among $k$ chains
with the $i$th chain receiving $a_i$ positions.
We already knew that the number of such linear extensions is the number of
complex solutions to a Wronski polynomial system for the chain polytope of
$P$. 

  Suppose now that $f(z)$ is a real polynomial with $r$ real roots and $c$ pairs
of complex conjugate roots, all distinct.
In every factorization of $f(z)$ into real polynomials, each conjugate pair of
roots must be distributed to the same polynomial.
This imposes stringent restrictions on the numbers of such real factorizations.

If every root of $f(z)$ is real, so that $c=0$, then the number of real
factorizations~\eqref{E:New_factor} is the multinomial
coefficient~\eqref{E:multinomialC}.
Also, there are no such factorizations if $f(z)$ has fewer than 
$|\{j\mid a_j\mbox{ is odd}\}|$ real roots.
In particular, the minimum number of real factorizations is 0 if more than one
$a_j$ is odd.
Recall that if $B\neq b_1+b_2+\dotsb +b_k$, then we have
\[
   \binom{B}{b_1,b_2,\dotsc,b_k}\ =\ 0\,.
\]

\begin{thm}\label{T:Fact_lower_bound}
 Suppose that $f(z)$ is a real polynomial of degree $a_1+\dotsb+a_k$ with distinct
 roots. 
 Let $n$ be the number of real factorizations~\eqref{E:New_factor} of $f$ 
 where $f_i$ has degree $a_i$.
 Then $n$ depends only on the number of real roots of $f(z)$ and satisfies

\[
   \binom{ \lfloor \frac{a_1+\dotsb+a_k}{2}\rfloor}%
   {\lfloor\frac{a_1}{2}\rfloor,\dotsc,\lfloor\frac{a_k}{2}\rfloor} 
    \ \leq \ n\ \leq \ 
    \binom{a_1+a_2+\dotsb+a_k}{a_1,a_2,\dotsc,a_k}\ .
\]
 The minimum is attained when $f(z)$ has at most one real root, and the maximum
 occurs when $f(z)$ has all roots real.
 Moreover, at most 
\[
   1 + \left\lfloor \frac{a_1}{2}\right\rfloor +
   \left\lfloor \frac{a_2}{2}\right\rfloor + \dotsb +
   \left\lfloor \frac{a_k}{2}\right\rfloor 
\]
 distinct values of $n$ can occur.
\end{thm}

For example, if $k=3$ and $(a_1,a_2,a_3)=(4,4,5)$, then $f(z)$ has degree 13.
The number $n$ of real factorizations of $f(z)$ into polynomials of degrees 4,4, and
5 as a function of the number of real roots $r$ of $f(z)$ is given in the table
below 
\begin{center}
  \begin{tabular}{|c||c|c|c|c|c|c|c|}\hline
   $r$ & 1 & 3 & 5 &  7 &  9 &  11 &   13\\\hline
   $n$ & 90&210&666&2226&7434&25410&90090\\ \hline
  \end{tabular}\vspace{5pt}
\end{center}

\noindent{\it Proof.}
  A factorization of a polynomial $f(z)$ with $r$ distinct real roots and $c$
distinct pairs of conjugate roots into real polynomials of degrees
$a_1,\dotsc,a_k$ is a distribution of the roots of $f$ among the factors where
the $i$th factor receives $a_i$ roots, and the conjugate pairs are
distributed to the same factor.

The upper bound was described previously, so we consider the lower bound.
The binomial coefficient lower bound vanishes when more than one $a_i$ is odd,
and we already observed that there are no real factorizations of $f$ in this
case.
If every $a_i$ is even and $f(z)$ has no real roots, then the root distribution
is enumerated by this binomial coefficient.
Lastly, if $a_i$ is the only odd number among $a_1,a_2,\dotsc,a_k$, and $f$ has
exactly one real root, that root must be given to the factor $f_i$.
If we replace $a_i$ by $a_i-1$ and this problem of distributing roots reduces to
the previous case. 

The last statement follows as $n$ depends only on the number of real
roots of $f(z)$ and $n=0$ unless $f(z)$ has at least 
$|\{j\mid a_j\mbox{ is odd}\}|$ real roots.
\hfill \raisebox{-3pt}{\epsfxsize=12pt\epsffile{figures/QED.eps}}\medskip

The number of real factorizations~\eqref{E:New_factor} is given by a generating
function. 
We thank Ira Gessel who explained this to us.

\begin{prop}
 The coefficient of $x_1^{a_1}\dotsb x_k^{a_k}$ in 
  $(x_1+\dotsb+x_k)^r(x_1^2 + \dotsb + x_k^2)^c$ is the number of
  factorizations 
\[
  f_1(z)\cdot f_2(z)\dotsb f_k(z)\ =\ f(z)
\]
 where $f(z)$ is real and has degree $r+2c=a_1+\dotsb+a_k$ with $r$
 distinct real roots and $c$ distinct pairs of complex conjugate roots, and 
 $f_i(z)$ is real and has degree $a_i$ for $i=1,\dotsc,k$.
\end{prop}

\noindent{\it Proof.}
 This is a standard use of generating functions, as described in Chapter 1
 of~\cite{St86a}. 
 We have $r$ red balls and $c$ cyan balls to distribute among $k$ boxes such
 that if $r_i$ is the number of red balls in box $i$ and $c_i$ is the number of
 cyan balls in box $i$, then $r_i+2c_i=a_i$.
 \hfill \raisebox{-3pt}{\epsfxsize=12pt\epsffile{figures/QED.eps}}\medskip

We relate the lower bound of Theorem~\ref{T:Fact_lower_bound} to
Theorem~\ref{T:Chain_Poly}. 

\begin{prop}\label{P:sign-imbalance_chain}
 Let $P$ be the incomparable union of chains of lengths 
 $a_1, a_2, \dotsc, a_k$.
 The  sign-imbalance of $P$ is
\[
  \sigma(a_1,a_2,\dotsc,a_k)\ :=\ 
  \binom{\lfloor\frac{a_i+\dotsb+a_k}{2}\rfloor}%
    {\lfloor\frac{a_1}{2}\rfloor,\dotsc,\lfloor\frac{a_k}{2}\rfloor}\,.
\]
 This equals zero unless at most one $a_i$ is odd.
\end{prop}

\noindent{\it Proof.}
 If we precompose a linear extension with the inverse of the extension where
 every element of the $i$th chain precedes every element of the ($i{+}1$)st
 chain, then we have identified the set of all linear extensions of $P$ with the
 set of minimal coset representatives $S^a$ of the subgroup 
 $S_{a_1}\times S_{a_2}\times\dotsb\times S_{a_k}$ of the symmetric group 
 $S_{a_1+\dotsb+a_m}$, which we call $(a_1,\dotsc,a_k)$-shuffles.
 The generating function for the distribution of lengths of these
 shuffles is the $q$-multinomial coefficient (the case $k=2$
 is~\cite[Prop. 1.3.7]{St86a}) 
\[
   \sum_{w\in S^a} q^{\ell(w)}\ =\ 
   \binom{a_1+a_2+\dotsb+a_k}{a_1,a_2,\dotsc,a_k}_q\ ,
\]
 where, if $k>2$, then 
 \begin{equation}\label{E:q-multinom}
   \binom{a_1+a_2+\dotsb+a_k}{a_1,a_2,\dotsc,a_k}_q\ =\ 
    \binom{a_1+a_2+\dotsb+a_{k-1}}{a_1,\dotsc,a_{k-1}}_q \cdot
    \binom{a_1+a_2+\dotsb+a_k}{a_1+\dotsb+a_{k-1},a_k}_q
 \end{equation}
 and $\binom{a+b}{a,b}_q$ is the $q$-binomial coefficient
 \begin{equation}\label{E:q-binom}
  \binom{a+b}{a,b}_q\ =\ 
   \frac{(1-q^{a+b})(1-q^{a+b-1})\dotsb(1-q^2)(1-q)}%
    {(1-q^a)\dotsb(1-q^2)(1-q)\cdot(1-q^b)\dotsb(1-q^2)(1-q)}\ .
 \end{equation}

 We evaluate the $q$-multinomial coefficient at $q=-1$ to compute the
 sign-imbalance of $P$.
 If $k$ is odd, then $1-q^k=2$ when $q=-1$.
 For even exponents, we have 
\[
  1-q^{2a}\ =\ (1-q^2)(1+q^2+q^4+\dotsb+q^{2a-2})
\]
 Now consider~\eqref{E:q-binom} when $q=-1$.
 If both $a$ and $b$ are odd, then~\eqref{E:q-binom} has one more factor with an
 even exponent in its numerator then in its denominator, and so it vanishes 
 when $q=-1$.
 Otherwise~\eqref{E:q-binom} has the same number of factors with even 
 exponents in its numerator as in its denominator, and so we cancel all
 factors of $(1-q^2)$.
 If we substitute $q=-1$, then each factor $(1-q^c)$ with odd exponent $c$
 becomes 2, and these cancel as there is the same number of such factors in the
 numerator and denominator.
 Since $(1+q^2+q^4+\dotsb+q^{2l-2})=l$ when $q=-1$, we see that
\[
  \binom{a+b}{a,b}_{q=-1}\ =\ \
  \binom{\lfloor\frac{a}{2}\rfloor+\lfloor\frac{b}{2}\rfloor}%
   {\lfloor\frac{a}{2}\rfloor,\lfloor\frac{b}{2}\rfloor}\ .
\]
 Applying~\eqref{E:q-multinom} to this formula completes the proof.
\hfill \raisebox{-3pt}{\epsfxsize=12pt\epsffile{figures/QED.eps}}\medskip

\begin{rem}
 By Theorem~\ref{T:Fact_lower_bound} and
 Proposition~\ref{P:sign-imbalance_chain}, the sign-imbalance of $P$ is the 
 sharp lower bound for the Wronski polynomial systems of the chain polytopes
 chain polytope of $P$.
 Thus Theorem~\ref{T:Chain_Poly} is sharp.
 Moreover, if the $a_i$ do not all have the same parity, then the hypotheses of 
 Theorem~\ref{T:Chain_Poly} do not hold, and in fact the toric variety $Y^+_P$
 is not orientable.
 Despite this, the conclusion of Theorem~\ref{T:Chain_Poly} does hold. 
\end{rem}

The ideas in Proposition~\ref{P:sign-imbalance_chain} can be used to compute the
sign-imbalance of a product of posets.
If $P$ is the incomparable union of posets $P_1,P_2,\dotsc,P_k$ with
$|P_i|=a_i$, then the linear extensions of $P$ are $(a_1,\dotsc,a_k)$-shuffles
of linear extensions of each component $P_i$.
If we let $\eta(P)$ be the number of linear extensions of a poset $P$, then 
we have the following corollary.

\begin{cor}\label{C:products}
 Let $P$ be as described.
 Then we have
\begin{eqnarray*}
   \eta(P)&=&\binom{a_1+a_2+\dotsb+a_k}{a_1,a_2,\dotsc,a_k}
   \cdot\prod_{i=1}^k \eta(P_i) \\
   \sigma(P)&=& \prod_{i=1}^k
   \sigma(P_i)\cdot\binom{\lfloor\frac{a_1+\dotsb+a_x}{2}\rfloor}%
          {\lfloor\frac{a_1}{2}\rfloor,\dotsc,\lfloor\frac{a_k}{2}\rfloor}\,.   
\end{eqnarray*}
\end{cor}

\begin{ex}
 The Grassmannian $G=G(2,2)$ has a sagbi degeneration to the toric variety
 $Y$ associated to the distributive lattice of order ideals on a product
 $C_2\times C_2$ of  two chains of length 2.
 Let $Z$ be the toric variety associated to the chain polytope of this
 poset.
 Since $C_2\times C_2$ is sign-balanced, the lower bound here is 0.

 If we take the product of $G$ with the projective plane, we obtain a
 variety to which Theorem~\ref{T:sagbi} applies.
 It has a sagbi degeneration into $Y\times\RP^2$, which
 is the toric variety of the distributive lattice of order ideals on the disjoint
 union of a chain $C_2$ of length 2 with $C_2\times C_2$.
 Similarly, the toric variety associated to the chain polytope of this
 poset is $Z\times\R\P^2$.
 By Corollary~\ref{C:products}, the Wronski polynomial systems on these
 varieties will have $30=2\cdot \binom{2+4}{2,4}$ complex solutions with at
 least 2 real. 

 The table below records the percentage that a given number of real roots was
 observed in Wronski polynomial systems on these varieties.
 The entries of 0 indicate values that were not observed.
\[
   \begin{tabular} {|c||c|c|c|c|c|c|c|c|c|c|c|c|c|c|c|c|}\hline
   \# real
   &{\bf 0}&{\bf 2}&{\bf 4}&{\bf 6}&{\bf 8}&{\bf 10}&{\bf 12}&{\bf 14}&
   {\bf 16}&{\bf 18}&{\bf 20}&{\bf 22}&{\bf 24}&{\bf 26}&{\bf 28}&{\bf 30}\\\hline
  $G\times \RP^2 $&0&11&55&24&5.3&1.6&1.2&1.4&0&0&0&0&0&0&0&.31\rule{0pt}{14pt}\\\hline
  $Y\times\RP^2$& 0& 2.2&25&14&23& 1.2&  .09&22&0&0&0&0.001& 0.001&.07& .01&12\rule{0pt}{14pt}\\\hline
  $Z\times\RP^2$&0& .07&6&33&4.6&1.3&2.9&39&0&0& .003& .01&1.5& .4& .37&10\rule{0pt}{14pt}\\\hline
\end{tabular}
\]
 We do not yet understand the apparent gaps in these data.
\end{ex}

\providecommand{\bysame}{\leavevmode\hbox to3em{\hrulefill}\thinspace}
\providecommand{\MR}{\relax\ifhmode\unskip\space\fi MR }
\providecommand{\MRhref}[2]{%
  \href{http://www.ams.org/mathscinet-getitem?mr=#1}{#2}
}
\providecommand{\href}[2]{#2}


\begin{thebibliography}{10}

\bibitem{Co95}
D.~Cox, \emph{The homogeneous coordinate ring of a toric variety}, J. Alg.
  Geom. \textbf{4} (1995), 17--50.

\bibitem{CEP82}
Corrado De~Concini, David Eisenbud, and Claudio Procesi, \emph{Hodge algebras},
  Ast\'erisque, vol.~91, Soci\'et\'e Math\'ematique de France, Paris, 1982,
  With a French summary. \MR{85d:13009}

\bibitem{DeKh00}
A.~I. Degtyarev and V.~M. Kharlamov, \emph{Topological properties of real
  algebraic varieties: {R}okhlin's way}, Uspekhi Mat. Nauk \textbf{55} (2000),
  no.~4(334), 129--212. \MR{1 786 731}

\bibitem{EG01a}
A.~Eremenko and A.~Gabrielov, \emph{The {W}ronski map and {G}rassmannians of
  real codimension 2 subspaces}, Comput. Methods Funct. Theory \textbf{1}
  (2001), no.~1, 1--25. \MR{2003h:26022}

\bibitem{EG01b}
\bysame, \emph{Degrees of real {W}ronski maps}, Discrete Comput. Geom.
  \textbf{28} (2002), no.~3, 331--347. \MR{2003g:14074}

\bibitem{Fu93}
William Fulton, \emph{Introduction to toric varieties}, Annals of Mathematics
  Studies, vol. 131, Princeton University Press, Princeton, NJ, 1993, The
  William H. Roever Lectures in Geometry. \MR{94g:14028}

\bibitem{Singular}
 G.-M. Greuel, G.~Pfister, and H.~Sch\"onemann, \emph{{\sc Singular} 2.0}, {A
  Computer Algebra System for Polynomial Computations}, Centre for Computer
  Algebra, University of Kaiserslautern, 2001, {\tt
  http://www.singular.uni-kl.de}.

\bibitem{Hibi}
T.~Hibi, \emph{Distributive lattices, affine semigroup rings and algebras with
  straightening laws}, Commutative Algebra and Combinatorics, Advanced Studies
  in Pure Mathematics, Vol.~11, North-Holland, 1987, pp.~93--109.

\bibitem{IKS03}
I.~Itenberg, V.~Kharlamov, and E.~Shustin, \emph{{Welschinger invariant and
  enumeration of real plane rational curves}}, International Mathematics
  Research Notes (2003), no.~49, 2639--2653.

\bibitem{J02}
Michael Joswig, \emph{Projectivities in simplicial complexes and colorings of
  simple polytopes}, Math. Z. \textbf{240} (2002), no.~2, 243--259.
  \MR{2003f:05047}

\bibitem{KM94}
M.~Kontsevich and {Yu}. Manin, \emph{{G}romov-{W}itten classes, quantum
  comomology, and enumerative geometry}, Comm. Math. Phys. \textbf{164} (1994),
  525--562.

\bibitem{Ko75}
A.G. Kushnirenko, \emph{A {N}ewton polyhedron and the number of solutions of a
  system of $k$ equations in $k$ unknowns}, Usp. Math. Nauk. \textbf{30}
  (1975), 266--267.

\bibitem{Kr1968}
L.~Kronecker, \emph{{L}eopold {K}ronecker's {W}erke}, Chelsea, NY, 1968.

\bibitem{Mi03}
Grigory Mikhalkin, \emph{Counting curves via lattice paths in polygons}, C. R.
  Math. Acad. Sci. Paris \textbf{336} (2003), no.~8, 629--634. \MR{1 988 122}

\bibitem{Mi04}
\bysame, \emph{Enumerative tropical algebraic geometry},  J. Amer. Math. Soc. 18
(2005), 313-377.  

\bibitem{So03}
Frank Sottile, \emph{Enumerative real algebraic geometry}, Algorithmic and
  quantitative real algebraic geometry (Piscataway, NJ, 2001), DIMACS Ser.
  Discrete Math. Theoret. Comput. Sci., vol.~60, Amer. Math. Soc., Providence,
  RI, 2003, pp.~139--179. \MR{1 995 019}

\bibitem{SS01}
Frank Sottile and Bernd Sturmfels, \emph{A sagbi basis for the quantum
  {G}rassmannian}, J. Pure Appl. Algebra \textbf{158} (2001), no.~2-3,
  347--366. \MR{2002a:13027}

\bibitem{St86b}
Richard Stanley, \emph{Two poset polytopes}, Discrete and Comput. Geom.
  \textbf{1} (1986), 9--23.

\bibitem{St86a}
\bysame, \emph{Enumerative combinatorics. {V}ol. 1}, Cambridge
  University Press, Cambridge, 1997, With a foreword by Gian-Carlo Rota,
  Corrected reprint of the 1986 original. \MR{98a:05001}

\bibitem{St02}
\bysame, \emph{{Some remarks on sign-balanced and maj-balanced
  posets}}, Advances in Applied Math., to appear.

\bibitem{Sturmfels_Alg}
Bernd Sturmfels, \emph{Algorithms in invariant theory}, Texts and Monographs in
  Symbolic Computation, Springer-Verlag, Vienna, 1993. \MR{94m:13004}

\bibitem{St94b}
\bysame, \emph{On the number of real roots of a sparse polynomial
  system}, Hamiltonian and gradient flows, algorithms and control, Fields Inst.
  Commun., vol.~3, American Mathematical Society, Providence, 1994,
  pp.~137--143.

\bibitem{Sturmfels_GBCP}
\bysame, \emph{Gr\"obner bases and convex polytopes}, American
  Mathematical Society, Providence, RI, 1996. \MR{97b:13034}

\bibitem{Wa96}
David~G. Wagner, \emph{Singularities of toric varieties associated with finite
  distributive lattices}, J. Algebraic Combin. \textbf{5} (1996), no.~2,
  149--165. \MR{97a:14057}

\bibitem{We03}
Jean-Yves Welschinger, \emph{Invariants of real rational symplectic 4-manifolds
  and lower bounds in real enumerative geometry}, C. R. Math. Acad. Sci. Paris
  \textbf{336} (2003), no.~4, 341--344. \MR{1 976 315}

\bibitem{Wh01}
Dennis~E. White, \emph{Sign-balanced posets}, J. Combin. Theory Ser. A
  \textbf{95} (2001), no.~1, 1--38. \MR{2002e:05151}

\end{thebibliography}
\end{document}